\def\Datum{December 8, 1998}

\magnification=\magstephalf
\hsize=16true cm\vsize=24.8true cm
\frenchspacing\parskip4pt plus 1pt

\newread\aux\immediate\openin\aux=\jobname.aux
\ifeof\aux\message{ <<< Run TeX a second time >>> }
\else\input\jobname.aux\fi\closein\aux
%\input \jobname.def

%   This is scr.def as of October 9, 1998
%
\parindent25pt
\newif\ifdraft\draftfalse
\newdimen\SkIp\SkIp=\parskip
\pageno=1
\long\def\fussnote#1#2{{\baselineskip=9pt
     \setbox\strutbox=\hbox{\vrule height 7 pt depth 2pt width 0pt}%
     \eightrm
     \footnote{#1}{#2}}}

\def\footnoterule{\kern-3pt
         \hrule width 2 true cm
         \kern 2.6pt}

\def\diagram{%
\def\normalbaselines{\baselineskip20pt\lineskip3pt\lineskiplimit3pt}
\matrix}
\def\mapright#1{\smash{\mathop{\hbox to 35pt{\rightarrowfill}}\limits^{#1}}}
\def\mapdown#1{\Big\downarrow\rlap{$\vcenter{\hbox{$\scriptstyle#1$}}$}}

 %neue Zeile ohne Indent und Skip

\font\gros=cmr17 scaled \magstep1
\font\gross=cmbx10 scaled\magstep1
\font\Gross=cmbx10 scaled\magstep2
\font\Sans=cmss10
\font\eightrm=cmr8
\font\eightsl=cmsl8
\font\eightsy=cmsy8
\font\eighti=cmmi8
\font\eightbf=cmbx8

\font\Icke=cmcsc10

\def\Bf{\noindent\bf}

%\font\tenmsa=msam10\font\sevenmsa=msam7\font\fivemsa=msam5
%\newfam\msafam
%\textfont\msafam=\tenmsa\scriptfont\msafam=\sevenmsa
%\scriptscriptfont\msafam=\fivemsa
%\def\msa{\fam\msafam}

\font\tenmsb=msbm10\font\sevenmsb=msbm7\font\fivemsb=msbm5
\newfam\msbfam
\textfont\msbfam=\tenmsb\scriptfont\msbfam=\sevenmsb
\scriptscriptfont\msbfam=\fivemsb
\def\msb{\fam\msbfam}

\def\rtimes{{\msb\char111}}

\font\tenmib=cmmib10 \font\eightmib=cmmib8 \font\fivemib=cmmib5
\newfam\mibfam
\textfont\mibfam=\tenmib \scriptfont\mibfam=\eightmib
\scriptscriptfont\mibfam=\fivemib
\def\mib{\fam\mibfam\tenmib}

\font\tenfr=eufm10 \font\eightfr=eufm8 \font\fivefr=eufm5
\newfam\frfam
\textfont\frfam=\tenfr \scriptfont\frfam=\eightfr
\scriptscriptfont\frfam=\fivefr
\def\frak{\fam\frfam\tenfr}

\font\tenbfr=eufb10 \font\eightbfr=eufb8 \font\fivebfr=eufb5
\newfam\bfrfam
\textfont\bfrfam=\tenbfr \scriptfont\bfrfam=\eightbfr
\scriptscriptfont\bfrfam=\fivebfr
\def\bfrak{\fam\bfrfam\tenbfr}

\def\Quot#1#2{\raise 1.5pt\hbox{$#1\mskip-1.2\thinmuskip$}\big/%
     \lower1.5pt\hbox{$\mskip-0.8\thinmuskip#2$}}

\def\klein{\eightrm\textfont1=\eighti\textfont0=\eightrm\def\sl{\eightsl}
\textfont2=\eightsy\baselineskip10.5pt\def\bf{\eightbf}}

\def\Kl#1{\raise1pt\hbox{$\scriptstyle($}#1\raise1pt\hbox{$\scriptstyle)$}}

\def\Bf{\noindent\bf}
\def\Exp{\mathop{\rm Exp}\nolimits}
\newcount\ite\ite=1
\def\0{\global\ite=1\1}
\def\1{\item{\rm(\romannumeral\the\ite)}\2}
\def\2{\global\advance\ite1}
\def\3#1{{\mib#1}}
\def\4#1{{\rm#1}}
\def\5#1{{\cal#1}}
\def\7#1{{\frak#1}}
\def\8#1{{\bfrak#1}}
\def\[#1{\mskip#1mu}
\def\]#1{\mskip-#1mu}

\def\To#1#2{(\romannumeral#1) $\Longrightarrow$ (\romannumeral#2)}

\def\endlist{\endlist}
\def\SperrRest{\afterassignment\SperrZeichen\let\next= }
\def\SperrZeichen{\ifx\next\endlist \let\next\relax
\kern-0.25em
\else \next \kern0.25em \let\next\SperrRest\fi\next}

\def\w{\ifmmode w^*\else $w^*$\fi}
\def\epsilon{\varepsilon}
\def\hat{\widehat}
\def\phi{\varphi}
\def\steil#1{\hbox{\rm\quad #1\quad}}
\def\Steil#1{\hbox{\rm\qquad #1\qquad}}
\def\p{p.\nobreak\hskip2pt}
\def\Im{\hbox{\rm Im}\[1}
\def\Re{\hbox{\rm Re}\[1}

\def\MOP#1{\expandafter\edef\csname #1\endcsname{%
   \mathop{\hbox{\Sans #1}}\nolimits}}

\MOP{Aff}
\MOP{Aut}
\MOP{Der}
\MOP{GL}
\MOP{SL}
\MOP{PSL}
\MOP{PSU}
\MOP{id}
\MOP{ad}
\MOP{Ad}
\MOP{U}
\MOP{O}
\MOP{SU}
\MOP{Spin}
\MOP{SO}
\MOP{Str}
\MOP{Tri}
\MOP{Sp}
\MOP{tr}
\MOP{Ker}
\MOP{Fix}
\def\AUT{\Aut_{\rm CR}}

\def\sqr#1#2{{\,\vcenter{\vbox{\hrule height.#2pt\hbox{\vrule width.#2pt
height#1pt \kern#1pt\vrule width.#2pt}\hrule height.#2pt}}\,}}

\def\qed{\hfill\ifmmode\sqr66\else$\sqr66$\par\fi\rm}
\def\qex{\hfill\ifmmode\sqr44\else$\sqr44$\par\fi\rm}

\def\One{{1\kern-3.8pt 1}}
\def\one{{1\kern-4.5pt 1\kern1pt}}

\def\CC{{\mathchoice
{\rm C\mkern-8mu\vrule height1.45ex depth-.05ex width.05em\mkern9mu\kern-.05em}
{\rm C\mkern-8mu\vrule height1.45ex depth-.05ex width.05em\mkern9mu\kern-.05em}
{\rm C\mkern-8mu\vrule height1ex depth-.07ex width.035em\mkern9mu\kern-.035em}
{\rm C\mkern-8mu\vrule height.65ex depth-.1ex width.025em\mkern8mu\kern-.025em}}}
\def\OO{{\mathchoice
{\rm O\mkern-11mu\vrule height1.42ex depth-.03ex width.05em\mkern9mu\kern-.05em}
{\rm O\mkern-9mu\vrule height1.45ex depth-.05ex width.05em\mkern9mu\kern-.05em}
{\rm O\mkern-9mu\vrule height1ex depth-.07ex width.035em\mkern9mu\kern-.035em}
{\rm O\mkern-9mu\vrule height.65ex depth-.1ex width.025em\mkern8mu\kern-.025em}}}
\def\RR{{\rm I\kern-1.6pt {\rm R}}}

\def\PP{{\rm I\!P}}
\def\BB{{\rm I\!B}}
\def\DD{{\rm I\!D}}
\def\KK{{\rm I\!K}}
\def\HH{{\rm I\kern-1.6pt {\rm H}}}
\def\TT{{\fam\msbfam\tenmsb T}}
\def\GG{{\fam\msbfam\tenmsb G}}
\def\ZZ{{\fam\msbfam\tenmsb Z}}

\newcount\nummer
\newcount\parno\parno=0

\def\Kap#1{{\par\bigbreak\global\advance\parno by1%
{\vskip3pt\noindent\gross\the\parno. #1\hfil}\vskip3pt}\rm\nummer=0\nobreak}

\def\KAP#1#2{\par\bigbreak\global\advance\parno by1%
\def\test{#1}\ifx\test\empty\else%
\expandafter\let\csname#1\endcsname=\relax%
\immediate\write\aux{\def\csname#1\endcsname{\the\parno}}%
\expandafter\xdef\csname#1\endcsname{\the\parno}\fi%
{\vskip3pt\noindent\gross\the\parno. #2\hfil}\vskip3pt\rm\nummer=0\nobreak}

\newwrite\aux

\def\Randmark#1{\vadjust{\vbox to 0pt{\vss\hbox to\hsize%
{\fiverm\hskip\hsize\hskip1em\raise 2.5pt\hbox{#1}\hss}}}}

\def\PrN{\the\parno.\the\nummer}

\def\Write#1#2{\global\advance\nummer1\def\test{#1}\ifx\test\empty\else%
\ifdraft\Randmark{#1}\fi\expandafter\let\csname#1\endcsname=\relax%
\immediate\write\aux{\def\csname#1\endcsname{\PrN}}%
\expandafter\xdef\csname#1\endcsname{\PrN}\fi#2}

\def\Num#1{\Write{#1}{\PrN}}
\def\Leqno#1{\Write{#1}{\leqno(\PrN)}}

\def\Joker#1#2{{\smallbreak\noindent\bf\Num{#1}\Icke~ #2.~}\parskip0pt\rm}
\def\SJoker#1#2{{\smallbreak\noindent\bf\Num{#1} #2.~}\parskip0pt\sl}
\def\Theorem#1{{\smallbreak\noindent\bf\Num{#1} Theorem.~}\parskip0pt\sl}
\def\Proposition#1{{\smallbreak\noindent\bf\Num{#1} Proposition.~}\parskip0pt\sl}
\def\Lemma#1{{\smallbreak\noindent\bf\Num{#1} Lemma.~}\parskip0pt\sl}
\def\Corollary#1{{\smallbreak\noindent\bf\Num{#1} Corollary.~}\parskip0pt\sl}

\def\Remark#1{{\smallbreak\noindent\bf\Num{#1} Remark.~}\parskip0pt}
\def\Example#1{{\smallbreak\noindent\bf\Num{#1} Example.~}}
\def\Definition#1{{\smallbreak\noindent\bf\Num{#1} Definition.~}}
\long\def\Proof{\smallskip\noindent\bf Proof.\parskip\SkIp\rm~~}
\def\Formend{\par\parskip\SkIp\rm}

\def\ruf#1{{{\expandafter\ifx\csname#1\endcsname\relax\xdef\flAG{}%
\message{*** #1 nicht definiert!! ***}\ifdraft\Randmark{#1??}\fi\else%
\xdef\flAG{1}\fi}\ifx\flAG\empty{\bf??}\else\rm\csname#1\endcsname\fi}}
\def\Ruf#1{{\rm(\ruf{#1})}}

\newcount\lit\lit=1
\def\Ref#1{\item{\the\lit.}\expandafter\ifx\csname#1ZZZ\endcsname\relax%
\message{ >>> \the\lit. = #1  <<< }\fi%
\expandafter\let\csname#1\endcsname=\relax%
\immediate\write\aux{\def\csname#1\endcsname{\the\lit}}\advance\lit1\ifdraft\Randmark{#1}\fi}

\def\Lit#1{\expandafter\gdef\csname#1ZZZ\endcsname{1}[\ruf{#1}]}

%\drafttrue
\ifdraft\footline={\hss\fiverm incomplete draft of \Datum\hss}
\else\nopagenumbers\fi

\headline={\ifnum\pageno>1\sevenrm\ifodd\pageno  Symmetric Cauchy-Riemann
manifolds\hss{\tenbf\folio}\else{\tenbf\folio}\hss{W. Kaup and D. Zaitsev
}\fi\else\hss\fi}

\immediate\openout\aux=\jobname.aux

\centerline{\Gross On symmetric Cauchy-Riemann manifolds}
\bigskip\bigskip
$$\vbox{\hsize=0.33\hsize
{\obeylines\everypar{\hfil}\parindent=0pt\parskip-3pt
\Icke Wilhelm Kaup
\vskip5pt
\sevenrm Mathematisches Institut
Universit\"at T\"ubingen
Auf der Morgenstelle 10
D-72076 T\"ubingen
Germany
 ~
kaup@uni-tuebingen.de
}}\qquad
\vbox{\hsize=0.33\hsize
{\obeylines\everypar{\hfil}\parindent=0pt\parskip-3pt
\Icke Dmitri Zaitsev
\vskip5pt
\sevenrm Mathematisches Institut
Universit\"at T\"ubingen
Auf der Morgenstelle 10
D-72076 T\"ubingen
Germany
 ~
dmitri.zaitsev@uni-tuebingen.de}}$$
\bigskip\bigskip

%{\narrower\klein\noindent  ABSTRACT.
%We study ..................\par}

\KAP{eins}{Introduction}

The Riemannian symmetric spaces
play an important role in different branches of mathematics.
By definition, a (connected) Riemannian manifold $M$ is called symmetric
if, to every $a\in M$, there exists an involutory isometric
diffeomorphism $s_a\colon M\to M$ having $a$ as isolated
fixed point in $M$ (or equivalently, if the differential
$d_a\!s_a$ is the negative identity on the the tangent space
$T_a=T_aM$ of $M$ at $a$). In case such a transformation $s_a$ 
exists for $a\in M$, it is uniquely determined and is the
geodesic reflection of $M$ about the point $a$. As a consequence,
for every Riemannian symmetric space $M$, the group $G=G_M$ 
generated by all symmetries $s_a$, $a\in M$, is a Lie group
acting transitively on $M$. In particular, $M$ can be identified with
the homogeneous space $G/K$ for some compact subgroup
$K\subset G$. Using the elaborate theory of Lie groups
and Lie algebras {\Icke E.$\,$Cartan} classified all Riemannian 
symmetric spaces.

The complex analogues of the Riemannian symmetric spaces
are the Hermitian symmetric spaces. By definition a Hermitian 
symmetric space is a Riemannian symmetric space $M$ together
with an almost complex structure on $M$ such that the
metric is Hermitian and such that every symmetry $s_a$
is holomorphic (i.e. satisfies the Cauchy-Riemann partial
differential equations with respect to the almost complex
structure). Also all Hermitian symmetric spaces were completely
classified by {\Icke E.$\,$Cartan}. In particular, every
Hermitian symmetric space $M$ can be written in a unique way
as an orthogonal direct product $M=M_+\times M_0\times M_-$
of Hermitian symmetric spaces $M_\epsilon$ with holomorphic
curvature of sign $\epsilon$ everywhere
(possibly of dimension 0, i.e. a single point).
$M_+$ is a compact simply connected complex manifold, $M_-$ is
a bounded domain in some $\CC^n$ and $M_0$ can be realized as
the flat space $\CC^m/\Omega$ for some discrete subgroup 
$\Omega\subset\CC^m$. In particular, the almost complex structure
of $M$ is integrable. Furthermore, there exists a remarkable
duality between symmetric spaces which for instance gives
a one-to-one correspondence between those of compact
type (i.e. $M=M_+$) and those of non-compact type (i.e. $M=M_-$),
see \Lit{HELG} for details.

A joint generalization of real smooth as well of complex manifolds
are the Cauchy-Riemann manifolds (CR-manifolds for short) or, more generally,
the CR-spaces, where the integrability condition is dropped and thus
also arbitrary almost complex manifolds are incorporated. These objects
are smooth manifolds $M$ such that at every point 
$a\in M$ the Cauchy-Riemann equations only apply in the direction
of a certain linear subspace $H_a\subset T_a$ of the tangent
space to $M$, see f.i. \Lit{BOGG} or section \ruf{zwei} for details.
The tangent space $T_a$ is an $\RR$-linear space while $H_a$, also
called the {\sl holomorphic} tangent space at $a\in M$, is a $\CC$-linear
space. The two extremal cases $H_a=0$ and $H_a=T_a$ for all $a\in M$
represent the two cases of smooth and of almost complex manifolds
respectively.

The main objective of this paper is to generalize the notion of
symmetry to the category of CR-spaces. It turns out that
for symmetries in this more general context
the requirement of isolated fixed points is no longer adequate. 
In fact, this would happen only for Levi-flat CR-spaces (see Proposition~\ruf{FLAT}) and hence would not be interesting.
Let us illustrate our concept on a simple example 
(compare also section \ruf{vier}). Consider $E\colon=\CC^n$, $n>1$,
with the standard inner product as a complex Hilbert space and
denote by $S\colon=\{z\in E:\|z\|=1\}$ the euclidean unit sphere
with Riemannian metric induced from $E$. Then $M$ is symmetric when
considered as Riemannian manifold. But $S$ also has a canonical
structure of a CR-manifold -- define for every $a\in S$ the
holomorphic tangent space $H_a$ to be the maximal complex
subspace of $E$ contained in $T_aS\subset E$, i.e. the
{\sl complex} orthogonal complement of $a$ in $E$. Then $T_a$ is the
orthogonal sum $H_a\oplus i\RR a$. It can be seen that for
every isometric CR-diffeomorphism $\phi$ of $S$, the differential 
$d_a\phi$ is the identity on $i\RR a$ as soon as it is the
negative identity on $H_a$, i.e. there does not exist a 
CR-symmetry of $S$ at $a$ in the strict sense. On the other hand,
the unitary reflection $s_a(z)\colon=2(a|z)-z$ defines an involutory 
isometric CR-diffeomorphism
of $S$ with differential at $a$ being the negative identity on $H_a$.
We call this the symmetry of the CR-manifold $S$ at $a$ 
and take it as a guideline for our general definition \ruf{SCR}.

Among all symmetric CR-manifolds we distinguish a large 
subclass generalizing the above example.
This class consists of the Shilov boundaries $S$ of bounded 
symmetric domains $D\subset\CC^n$
in their circular convex realizations (Theorem~\ruf{MINI}).
A remarkable feature of these CR-submanifolds is 
the fact that various geometric and analytic constructions,
hard to calculate in general, can be obtained here in very explicit forms.
We illustrate this on the case of polynomial and rational convex hulls.

Recall that the polynomial (resp. rational) convex hull of a compact subset $K\subset\CC^n$
is the set of all $z\in\CC^n$ such that $|f(z)|\le \sup_K |f|$ for every polynomial $f$
(resp. every rational function $f$ holomorphic on $K$). 
If $K$ is a connected real-analytic curve, $p(K)\backslash K$ is a complex analytic subset of $\CC^n\backslash K$,
due to {\Icke J.$\,$Wermer} \Lit{WERA}, where $p(K)$ denotes the polynomial convex hull. 
Later, the analyticity of $p(K)\backslash K$ 
was proved by {\Icke H.$\,$Alexander} \Lit{ALO}, \Lit{ALN} for compact sets of finite length 
and recently by {\Icke T.C.$\,$Dinh} \Lit{DINH} for rectifiable closed (1,1)-currents under very weak assumptions
(see also {\Icke E.$\,$Bishop} \Lit{BISH}, {\Icke G.$\,$Stolzenberg} \Lit{STO} 
and {\Icke M.G.$\,$Lawrence}\Lit{LAW} for related results). 
On the other hand, if $K$ is not a smooth submanifold, $p(K)\backslash K$ is not analytic in general
(see e.g. {\Icke G.$\,$Stolzenberg} \Lit{STOE} or {\Icke J.$\,$Wermer} \Lit{WERE}). 

In the present paper we calculate the polynomial and the rational convex hulls of $S$, 
where $S\subset \partial D$ is a Shilov boundary as above (see Corollary~\ruf{HULLS}).
Here it turns out that $p(S)\backslash S$ is not necessarily analytic but rather a 
{\sl submanifold with ``real-analytic corners''},
even though $S$ itself is a connected real-analytic submanifold.
Similar is the behaviour of the rational convex hull of $S$.
Both hulls are canonically stratified into real-analytic CR-submanifolds 
such that the (unique) stratum of the highest dimension is complex for the polynomial and Levi-flat for
the rational convex hull.

The paper is organized as follows.
Preliminaries on CR-spaces are given in section \ruf{zwei}.
In section \ruf{drei} we introduce symmetric CR-spaces 
and establish their main properties:
{\sl The uniqueness of symmetries and the transitivity of the spanned group}.
Example \ruf{EXAP} (a generalized Heisenberg group)
shows that there exist symmetric CR-manifolds $M$
of arbitrary CR-dimension and arbitrary 
CR-codimension having arbitrary Levi form at a given point.

In section \ruf{vier} we study more intensively the unit 
sphere $S$ in the complex space $\CC^n$,
some symmetric domains in $S$ and their coverings.
In particular, we obtain uncountable families of 
pairwise non-isomorphic symmetric CR-manifolds
(see Example \ruf{EXAX}). In section \ruf{fuenf} 
we associate to every symmetric CR-space $M$
a canonical fibration and discuss the situation 
when the base is a symmetric CR-space itself.
In fact, every symmetric CR-space can be obtained in this way.
As mentioned above, the underlined CR-structure of 
a symmetric CR-space does not need to be integrable.
In section \ruf{fuenfa} we give a construction 
principle for symmetric CR-spaces in terms of
Lie groups and illustrate this with various examples.
In section \ruf{sechs} we give Lie theoretic 
conditions for $M$ to be embeddable into
a complex manifold. We show by an example 
(see \ruf{VORH}) that this in general is
not possible -- in contrast to the case 
of Hermitian symmetric spaces. 

Finally, in section \ruf{sieben}, we consider
symmetric CR-manifolds
arising from bounded symmetric domains $D\subset\CC^n$.
To be a little more specific, we assume without loss
of generality that $D$ is realized as bounded circular
convex domain in $\CC^n$. Then it is known that the
Shilov boundary $S\subset \partial D$ of $D$ (which coincides here with
the set of all extreme points of the convex set $\overline D$)
is an orbit of the group $\Aut(D)$ of all biholomorphic 
automorphisms of $D$. Furthermore, every maximal compact subgroup $K$
of $\Aut(D)$ still acts transitively on $S$, and, with respect to a
suitable $K$-invariant Hermitian metric, $S$ is a symmetric
CR-manifold. Our main result states: In case $D$ does not have a
factor of tube type, (i) every smooth CR-function on $S$ extends
to a continuous function on $\overline D$ holomorphic on $D$, and
(ii) the group $\AUT(S)$ of all smooth CR-diffeomorphisms of $S$
coincides with the group $\Aut(D)$.

{\parindent0pt\noindent
{\bf Notation.}
For every vector space $E$ over the base field $\KK=\RR$ or $\KK=\CC$
we denote by $\5L(E)$ the space of all $\KK$-linear endomorphisms of $E$
and by $\GL(E)\subset\5L(E)$ the subgroup of all invertible operators.
More generally, for every total subset $S\subset E$ put
$\GL(S)\colon=\{g\in\GL(E):g(S)=S\}$
and denote by $\Aff(S)$ the group of all affine transformations
of $E$ mapping $S$ onto itself.

$\KK^{n\times m}$ is the $\KK$-Hilbert space of all $n\times m$-matrices
over $\KK$ ($n=$ row-index) with the inner product $(u|v)=\tr(u^*v)$ and
$u^*=\overline u'\in\KK^{m\times n}$ the adjoint.

By a complex structure we always understand a linear operator $J$
on a real vector space with $J^2=-\id$. If misunderstanding
is unlikely we simply write $ix$ instead of $J(x)$.

For complex vector spaces $U,V,W$ a sesqui-linear map $\Phi\colon
U\times V\to W$ is always understood to be conjugate linear in the
first and complex linear in the second variable.

With $\U(n)$, $\O(n)$, $\Sp(n)$ etc. we
denote the unitary, orthogonal and symplectic
groups (see \Lit{HELG} for related groups). In particular, we put
$\TT:=\U(1)=\exp(i\RR)$ and $\RR^+\colon=\exp(\RR)$.
For every topological group $G$ we denote by $G^0$ the connected identity
component. A continuous action of $G$ on a
locally compact space $M$ is called proper
if the mapping $G\times M\to M\times M$ defined by $(g,a)\mapsto (g\Kl a ,a)$ is
proper, i.e. pre-images of compact sets are compact.

For every set $S$ and every map $\sigma\colon S\to S$ we denote
by $\Fix(\sigma)\colon=\{s\in S:\sigma\Kl s=s\}$ the set of all
fixed points.

For Lie groups $G$, $H$ etc., the corresponding Lie algebras are denoted by small
Gothic letters $\7g$, $\7h$ etc. For linear subspaces $\7m,\7n\subset\7g$ we
denote by $[\7m,\7n]$ always the linear
span of all $[x,y]$ with $x\in\7m$, $y\in\7n$.

}

\KAP{zwei}{Preliminaries}

Suppose that $M$ is a connected smooth manifold of (finite real) dimension
$n$. Denote by $T_a\colon=T_aM$, $a\in M$, the tangent space which
is a real vector space of dimension $n$. An {\sl almost Cauchy-Riemann structure}
(almost CR-structure for short) on $M$ assigns
to every $a\in M$ a linear subspace $H_a=H_aM\subset T_a$
(called the {\sl holomorphic tangent space} to $M$ at $a$)
together with a complex structure
on $H_a$ in such a way that $H_a$ and the complex structure depend
smoothly on $a$. Smooth dependence can be expressed in the following way:
{\sl Every point of $M$ admits an open neighbourhood $U\subset M$ together
with a linear endomorphism $j_a$ of $T_a$ for every $a\in U$ such
that $-j_a^2$ is a projection onto $H_a$ with $j_av=iv$ for all $v\in H_a$,
and $j_a$ depending smoothly on $a\in U$}. Then, in particular, all $H_a$
have the same dimension. A connected smooth manifold together with an
almost CR-structure on it is called in the sequel an {\sl almost CR-manifold}, or
a {\sl CR-space} for short. For more details on (almost) CR-manifolds see e.g.
\Lit{BER}, \Lit{BOGG}, \Lit{CHIR}, \Lit{JACO}, \Lit{TUCR}.

In the following $M$ always denotes a CR-space.
Denote by $\7V=\7V(M)$ the Lie algebra
of all smooth vector fields on $M$ and by $\7H=\7H(M)$ the
subspace of all vector fields $X$ with $X_a\in H_a$ for
all $a\in M$. Then, for all $a\in M$, $\7H_a\colon=\{X_a:X\in\7H\}=H_a$ holds
and $(JX)_a=i(X_a)$ canonically defines a complex structure $J$ on $\7H$.
Define inductively $$\7H^k\colon
=\7H^{k-1}+[\7H,\7H^{k-1}]\,,\steil{where}\7H^1\colon=\7H
\steil{and}\7H^j\colon=0\;\;\hbox{\rm for}\;\;j\le0\,.\Leqno{HKKK}$$
Then $[\7H^r,\7H^s]\subset\7H^{r+s}$ holds for all integers $r,s$
and $\7H^\infty\colon=\bigcup_k\7H^k$ is the Lie
subalgebra of $\7V$ generated by $\7H$.
Also, we call the quotient vector spaces
$$T_a^r\colon=T_a/H_a\steil{and}T_a^{rr}\colon=T_a/\7H^\infty_a$$
the {\sl real} and the {\sl totally real} part of $T_a$.
The complex dimension of $H_a$ and the real dimension of
$T_a^r$ do not depend on $a\in M$ - they are called the {\sl
CR-dimension} and the {\sl CR-codimension} of $M$.
Finally, $M$ is called {\sl minimal} (in the sense of {\Icke Tumanov} \Lit{TUMA})
if $U=N$ holds for every domain
$U\subset M$ and every closed smooth submanifold
$N\subset U$ with $H_aM\subset T_aN$ for all $a\in N$.
It is known that in case $M$ is real-analytic (all
CR-spaces we discuss later will
have this property, see \ruf{REAN}) minimality
is equivalent to $T_a^{rr}=0$ for all $a\in M$, i.e. to the finite type
in the sense of {\Icke Bloom-Graham} \Lit{BG}. $M$ is called {\sl totally real} if $M$ has CR-dimension 0. 
The CR-spaces of CR-codimension 0, i.e. satisfying $H_aM=T_aM$, are also
called almost complex manifolds.

The CR-spaces form a category in a natural way. By definition,
a smooth map $\phi\colon M\to N$ of CR-spaces is called a CR-{\sl map} if,
for every $a\in M$ and $b=\phi(a)\in N$,
the differential $d_a\]1\phi\colon T_aM\to T_bN$ maps $H_aM$
complex linearly into $H_bN$. For every $M$ we denote
by $\Aut(M)$ or $\AUT(M)$ the group of all CR-diffeomorphisms
$\phi\colon M\to M$ and endow this group with the compact open
topology.

Suppose, $N$ is a CR-space and $M\subset N$ is a submanifold.
We call $M$ a CR-{\sl subspace} of $N$ if the dimension of
$H_aM\colon=(T_aM\cap H_aN)\,\cap\,i(T_aM\cap H_aN)$ does not depend
on $a\in M$. Then $M$ is a CR-space with the induced CR-structure.
A CR-space $M$ is called {\sl integrable} or a {\sl CR-manifold}
if the following integrability condition is satisfied:
$$Z\colon=[JX,Y]+[X,JY]\in\7H\steil{and}
JZ=[JX,JY]-[X,Y]\steil{for all}X,Y\in\7H\,.\Leqno{INTE}$$
In the special case, where $M$ is real-analytic (which includes that
the holomorphic tangent space $H_aM$ depends in a real-analytic
way on $a\in M$) it is known (compare f.i. \Lit{BER} or \Lit{BOGG}) that \Ruf{INTE}
is equivalent to the existence of local realizations of $M$
as a CR-submanifold of some $\CC^n$.
In that case there even exist a complex manifold $N$
and a (global) realization of $M$ as a real-analytic CR-submanifold
of $N$ which is generic, i.e. $T_aM+iT_aM=T_aN$ for all $a\in M$ (see \Lit{ANDR}).

An important invariant of a CR-space $M$ is the so-called {\sl Levi form}
defined at every point $a\in M$ in the following way. Denote by 
$\pi_a\colon T_a\to T_a/H_a$ the canonical projection. 
Then it is easy to see that there exists a map
$$\omega_a\colon H_a\times H_a\;\to\;T_a/H_a\steil{with}\omega_a(X_a,Y_a)=\pi_a([X,Y]_a)
\steil{for all}X,Y\in\7H$$
(in the proof of Proposition \ruf{UUNI} a more general satement actually
will be shown). Then $\omega_a$ is $\RR$-bilinear and skew-symmetric.
For all $\epsilon,\mu=\pm1$, there exist uniquely determined $\RR$-bilinear maps 
$$\omega_a^{\epsilon\mu}\colon H_a\times H_a\to (T_a/H_a)\otimes_\RR\CC\steil{with}
\omega_a^{\epsilon\mu}(sx,ty)=s^\epsilon t^\mu \omega_a^{\epsilon\mu}
(x,y)$$
for all $s,t\in\TT$, $x,y\in H_a$
such that $\omega_a=\sum \omega_a^{\epsilon\mu}$.
\Definition{LEVI} The sesqui-linear part $L_a\colon=\omega_a^{-1,1}$ of $\omega_a$
is called the {\sl Levi form} of $M$ at $a$, i.e.
$$4\,L_a(x,y)=\big(\omega_a(x,y)+\omega_a(ix,iy)\big)\,+\,i\big(\omega_a(ix,y)-\omega_a(x,iy)\big)$$
and in particular $2L_a(x,x)=i\omega_a(ix,x)\in i(T_a/H_a)$ for all $x,y\in H_a$.
The convex hull of $\{L_a(x,x):x\in H_a\}$ is called the {\sl Levi cone} at $a\in M$ and its interior always refers to the linear space $i(T_a/H_a)$. The CR-space
$M$ is called {\sl Levi flat} if $L_a=0$ holds for every $a\in M$.
\Formend
\noindent
Denote by  $^*$ the conjugation of $(T_a/H_a)\otimes_\RR\CC$
given by $(\xi+i\eta)^*\colon=(-\xi+i\eta)$ for all $\xi,\eta\in T_a/H_a$.
The following statement is obvious.
\Lemma{} The Levi form $L_a$ is $^*$-Hermitian,
that is, $L_a(x,y)=L_a(y,x)^*$ for all $x,y\in H_a$.
In the case, the  integrability condition \Ruf{INTE} holds,
we have $\omega_a^{1,1}=\omega_a^{-1,-1}=0$ and 
$$2\,L_a(x,y)\;=\;\omega_a(x,y)+ i\omega_a(ix,y)\,.\eqno{\qex}$$
\Formend
Various authors call $-2iL_a$ the Levi form.
In \Lit{BOGG} $L=L_a$ is called the {\sl extrinsic Levi form}.
If $M$ is a CR-subspace of a complex manifold $U$, $(T_a/H_a)\otimes_\RR\CC$
can be canonically identified with a complex subspace of $T_aU/H_aM$.
Then, for every $x\in H_a$ the vector $L(x,x)\in iT_a/H_a$ is transversal to 
$M$ and
points into the `pseudoconvex direction' of $M$.
{\Icke Bogges} and {\Icke Polking} \Lit{BP} proved that all 
CR-functions on $M$ extend holomorphically to a wedge in $U$
`in the direction of the Levi form'. This was generalized by
{\Icke Tumanov} \Lit{TUMA} to the case $M$ is minimal,
whereas {\Icke Baouendi} and {\Icke Rothschild} \Lit{BRMIN} proved the necessity of
the minimality condition.

Let us illustrate the Levi form at a simple example: 
Let $M\colon=\{(z,w)\in\CC^2:z\overline z+w\overline w=1\}$
be the euclidean sphere and put $a\colon=(1,0)$, $e\colon=(0,1)$.
Then $M$ is a CR-submanifold with $T_aM=i\RR a\oplus\CC e$ and
$H_aM=\CC e$. We identify $(T_aM/H_aM)\otimes_\RR\CC$
in the obvious way with the complex line $\CC a\subset\CC^2$. Consider the
vector field $X\in\7H$ defined by $X_{(z,w)}=(-\overline w,\overline z)$
for all $(z,w)\in M$.
Then $[JX,X]_a=2ia=\omega_a(ie,e)$ and hence $L_a(e,e)=-a$. This vector
points from $a\in M$ into the interior of the sphere $M$.

For the rest of the section assume that on $M$ there is given
a {\sl Riemannian metric}, i.e. every tangent space $T_a$
is a real Hilbert space with respect to an inner product
$\langle u|v\rangle_a$ depending smoothly on  $a\in M$.
Then we call $M$ an {\sl Hermitian CR-space}
if the metric is compatible with the CR-structure in the sense that
$\|iv\|=\|v\|$ holds for all $a\in M$, $v\in H_a$, where 
$\|v\|=\sqrt{\langle v,v\rangle}$. In particular, every $H_a$ is a
complex Hilbert space. If the CR-structure of a Hermitian CR-space $M$ is integrable,
we call $M$ an {\sl Hermitian CR-manifold}.
In the case of vanishing CR-codimension, i.e. $H_aM=T_aM$,
Hermitian CR-spaces are usually called {\sl Hermitian manifolds}.

For every $a\in M$ and every integer $k\ge0$,
let $\7H^k$ be as in \Ruf{HKKK} and let $H_a^k\subset T_a$ 
be the orthogonal complement of $\7H_a^{k-1}$
in $\7H_a^k$. Then $H_a^0=0$ and $H_a^1=H_a$ is the holomorphic
tangent space in $a$. If we denote by $H_a^{-1}\cong T_a^{rr}$ the
orthogonal complement of $\7H^\infty_a$ in $T_a$ we obtain the following orthogonal decomposition
$$T_a=\bigoplus_{k\ge-1}H_a^k\;.\Leqno{DECO}$$
Denote by $\pi_a^k\in\5L(T_a)$ the orthogonal projection onto $H_a^k$.
Then $a\mapsto\pi_a^k$ defines
a (not necessarily continuous) tensor field $\pi^k$ of type $(1,1)$ on $M$.
For later use we define
$$T_a^+\colon=\bigoplus_{k{\rm~even}}H_a^k\steil{and}
T_a^-\colon=\bigoplus_{k{\rm~odd}}H_a^k\;.\Leqno{DEOR}$$

The Hermitian CR-spaces form a category together with the
contractive CR-mappings as morphisms (i.e. $\|d_a\phi\Kl v \|\le\|v\|$
for all $a\in M$, $v\in T_aM$).
We always denote by $I_M\subset\AUT(M)$ the closed subgroup of
all isometric CR-diffeomorphisms. Then it is known that
$I_M$ is a Lie group acting smoothly and properly on $M$. 
In particular, $I_M$ has
dimension $\,\le n(n\]4+\]42)+m(m\]4+\]41)/2$, where $M$ has
CR-dimension $n$ and CR-codimension $m$. The full
CR-automorphism group $\AUT(M)$ can be infinite-dimensional.
In case of vanishing CR-dimension
we have the full sub-category of (connected) Riemannian manifolds and in case
of vanishing CR-codimension we have the full sub-category of
Hermitian manifolds. In both sub-categories there exists the classical
notion of a symmetric space. In the following we want to extend this
concept to arbitrary Hermitian CR-spaces.

\KAP{drei}{Symmetric CR-spaces.}

\Definition{STRY} Let $M$ be an Hermitian CR-space and let
$\sigma\colon M\to M$ be an isometric CR-diffeomorphism.
Then $\sigma$ is called a
{\sl symmetry} at the point $a\in M$ (and $a$ is called a {\sl symmetry point} of $M$)
if $a$ is a (not necessarily isolated) fixed point of $\sigma$  and if the differential
of $\sigma$ at $a$ coincides with the negative identity on the subspace
$H_a^{-1}\!\oplus H_a^1$ of $T_a$.
\Formend

\Proposition{UNIQ} At every  point of $M$ there exists at most one symmetry.
Furthermore, every symmetry is involutive.
\Formend\noindent The statement is an easy consequence of the following.
\SJoker{UUNI}{Uniqueness Theorem} Let $\phi,\psi$ be isometric CR-diffeomorphisms of
the Hermitian CR-space $M$ with $\phi\Kl a=\psi\Kl a $ for
some $a\in M$. Then $\phi=\psi$ holds if the differentials
$d_a\phi$ and $d_a\psi$ coincide on
the subspace $H_a^{-1}\!\oplus H_a^1$ of $T_a=T_aM$.
\Proof Without loss of generality we may assume that
$\psi$ is the identity transformation of $M$.
Then $a$ is a fixed point of $\phi$ and by a well known
fact (compare f.i. \Lit{HELG} \p62)
we only have to show that the differential
$\lambda\colon=d_a\phi$ is the identity on $T_a$.
Now $\lambda(X_a)=(\tau X)_a$ holds for every vector field $X\in\7V$,
where $\tau$ is the Lie automorphism of $\7V$ induced by
$\phi$. For all $r,s>0$ and $k\colon=r+s$, every smooth function
$f$ on $M$ and all $X\in\7H^{r}$, $Y\in\7H^{s}$ the formula
$$\pi_a^k\big([fX,Y]_a\big)\;=\;f\Kl a \cdot\pi_a^k\big([X,Y]_a\big)$$
is easily derived, where the orthogonal projection $\pi_a^k$ is defined
as above. On the other hand, every $X\in\7H^r$
with $X_a=0$ can be written as finite sum $X=X_0+f_1X_1+\cdots+f_mX_m$
with $X_0\in\7H^{r-1}$,
$X_1,\dots,X_m\in\7H^r$ and smooth functions $f_1,\dots,f_m$
on $M$ vanishing in $a$. Therefore, $\pi_a^k\big([X,Y]_a\big)$
 only depends on the
vectors $\pi_a^r(X_a)$ and $\pi_a^s(Y_a)$ for $X\in\7H^{r}$, $Y\in\7H^{s}$.
\hfill\break
Since $\lambda$ is
an isometry of $T_a$ and $\tau$ leaves invariant all
subspaces $\7H^k\subset\7V$ also all $H_a^k$ must be left
invariant by $\lambda$. We show by induction on $k$
that actually $\lambda$ is the identity on every $H_a^k$.
For $k\le1$ this follows from the assumptions.
For $k>1$, fix $X\in\7H$, $Y\in\7H^{k-1}$ and consider
the vector $v\colon=\pi_a^k\big([X,Y]_a\big)\in H_a^k$.
By induction hypothesis we then have $(\tau X)_a=X_a$, $(\tau Y)_a=Y_a$
and hence
$$\lambda\Kl v=\pi_a^k(\lambda\Kl v )=\pi_a^k\big((\tau[X,Y])_a\big)=\pi_a^k\big([\tau X,
\tau Y]_a\big)=\pi_a^k\big([X,Y]_a\big)=v\,.\eqno{\qed}$$
\Remark{SYHA} The proof of Proposition \ruf{UUNI} shows that for every
symmetry $\sigma$ of $M$ at $a$ the differential $\lambda=d_a\sigma$
satisfies $\lambda\Kl v=(-1)^kv$ for every $v\in H_a^k$ and every
$k\ge-1$, i.e. $\lambda=\Sigma_k(-1)^k\pi_a^k$,
or equivalently, $\,\Fix(\pm\lambda)=T_a^\pm$.
\qex

\Definition{SCR} A connected Hermitian CR-space $M$ is called
{\sl symmetric} (or an SCR-space for short)
if every $a\in M$ is a symmetry point.
The corresponding symmetry at $a$ is denoted by $s_a$.
\qex

In the sequel we adopt the following notation: For a given
SCR-space $M$ we denote as in section \ruf{zwei} by $I=I_M$ the
Lie group of all isometric CR-diffeomorphisms of $M$.
Let $G=G_M$ be the closed subgroup of $I$ generated by all symmetries $s_a$, $a\in M$.
Fix a base point $o\in M$ and denote by
$K\colon=\{g\in G:g\Kl o=o\}$ the isotropy subgroup at $o$.

\Proposition{SYMM} $G$ is a Lie group
acting transitively and properly on $M$. The connected
identity component $G^0$ of $G$ has index $\le2$ in $G$
and coincides with the closed subgroup of $I_M$ generated by all
transformations $\, s_a\circ s_b$ with $a,b\in M$. The isotropy
subgroup $K$ is compact and $M$ can be
canonically identified with the homogeneous manifold
$G/K$ via $g\Kl o \mapsto gK$.
$M$ is compact if and only if $G$ is a compact Lie group.
\Proof There exists an open subset $U\ne\emptyset$ of $M$
such that for every $k\ge-1$ the dimension of $H_a^k$
does not depend on $a\in U$. Therefore, every tensor field
$\pi^k$ is smooth over $U$, i.e. the orthogonal
decomposition \Ruf{DECO} depends smoothly on $a$ as long
as $a$ stays in $U$. We may assume without loss of generality
that for a fixed $\epsilon>0$ and every $a\in U$, the exponential mapping
$\Exp_a$ is defined on the open ball $B_a$ of
radius $\epsilon$ about the origin in $T_a$
and that $\Exp_a$ is a diffeomorphism from $B_a$
onto a neighbourhood $N_a\subset M$ of $a$. Every
isometric diffeomorphism $\phi\in I_M$ is linear in local normal
coordinates, more precisely, for every $a\in U$ with $c\colon=\phi\Kl a \in U$,
the diagram
$$\diagram{B_a&\mapright{\displaystyle d_a\phi} &B_c\cr
\mapdown{\displaystyle\hskip-29pt\Exp_a} & &\mapdown{\[4\displaystyle\Exp_c} \cr
N_a&\mapright{\displaystyle\phi} &N_c\cr}$$
commutes. Since $d_a\[2s_a=\sum_k(-1)^k\pi_a^k$ depends smoothly on $a\in U$
and since $G$ consists of isometries 
this implies the smoothness of the mapping $U\times U\to M$ defined by $(a,b)\mapsto s_a(b)$.
Now fix $u\in U$ and denote by $A\colon=G\Kl u $ the orbit of $u$
under the group $G$. Then $A$ is a closed smooth submanifold of $M$
since $G$ acts properly on $M$. Fix $a\in A$
and $v\in T_a^-M$ arbitrarily. Choose a smooth curve $\gamma\colon[0,1]\to U$
with $\gamma(0)=a$ and $\gamma'(0)=v$.
Then $\alpha(t)=s_{\gamma(t)}\Kl a $ defines a smooth curve in $A$
with $\alpha(0)=a$ and $\alpha'(0)=2v$.
This proves $T_a^-M\subset T_aA$ and hence $T_c^{-1}M\subset T_cA$ for all $c\in A$
since $a\in A$ was arbitrarily chosen. Now fix a vector $w\in T_a^+M$.
Then there exist smooth vector fields $X^1,\dots,X^{2k}$
on $M$ such that $X_c^j\in T_c^-M$ for all
$1\le j\le2k$, $c\in A$ and such that $w=\sum_{j=1}^k[X^j,X^{k+j}]_a$.
But we know already that
all $X^j$'s are tangent to the submanifold $A\subset M$, i.e. all their brackets are
tangent to $A$ and therefore $w\in T_aA$. This implies $T_aA=T_aM$ and $A=M$ since
$A$ is closed in $M$. Therefore $G$ acts transitively on $M$.
\qed

The proof of Proposition \ruf{SYMM} shows that an Hermitian CR-space
$M$ is already symmetric as soon as the set of symmetry points of $M$
has an interior point in $M$. The transitivity of the $G$-action
has several consequences.

\Corollary{REAN} Every SCR-space $M$ has a unique structure
of a real-analytic CR-manifold in such a way that every isometric diffeomorphism
of $M$ is real-analytic. All tensors $\pi^k$ are real-analytic on $M$
and also the mapping $a\mapsto s_a$ from $M$ to $G$ is real-analytic.
In particular, the dimension of $H_a^k\subset T_a$
does not depend on $a\in M$ for every $k$.
\qex

In the following we will always consider SCR-spaces
as real-analytic manifolds according to \ruf{REAN}.
The number
$\kappa=\kappa(M)\colon=\max\{k\ge-1:H_a^k\ne0\}$ does not depend on $a\in M$.
Example \ruf{SYIN} will show that arbitrary values of $\kappa\ne0$ occur.
$\CC^n\!\!\times\!\RR^m$ is a Levi flat CR-submanifold of $\CC^{n+m}$ and
as another corollary of \ruf{SYMM} we have:
\Proposition{FLAT} Let $M$ be a symmetric CR-space with CR-dimension $n$
and CR-codimension $m$. Then for every $a\in M$ the following
conditions are equivalent.
\0 $a$ is an isolated fixed point of the the symmetry $s_a$.
\1 $M$ is Levi flat.
\1 $M$ is locally CR-isomorphic to an open subset of $\CC^n\!\!\times\!\RR^m$.
In particular, $M$ is a CR-manifold.
\Proof \To12. The differential $d_as_a$ is the identity on the subspace
$H_a^2\subset T_a$ due to Remark \ruf{SYHA}. Therefore, if $a$ is isolated in $\Fix(s_a)$, we have
$H_a^2=0$ and hence $L_a=0$. By homogeneity then the Levi form vanishes
at every point of $M$, i.e. $M$ is Levi flat.\hfill\break
\To23. Suppose that $M$ is Levi flat. Then the holomorphic tangent
spaces form an involutive distribution on $M$ and define a foliation
of $M$. Let $N$ be the leaf through $a$, i.e. the maximal connected
immersed smooth submanifold $N$ of $M$ with $T_cN=H_cM$ for all $c\in N$.
Then $N$ is an Hermitian almost complex manifold in the leaf topology
invariant under every symmetry $s_c$, $c\in N$. Therefore, $N$ is
an Hermitian symmetric space and in particular a complex manifold,
see \Lit{HELG}. But $M$ locally is CR-isomorphic to a direct
product $U\times V$, where $U\subset N$ and $V\subset\RR^m$ are
open subsets.\hfill\break
\To31. Condition (iii) implies $H_a^{-1}\oplus H_a^1=T_a$ and hence
$d_as_a=-\id$, i.e. $a$ is an isolated fixed point of $s_a$.
\qed

Every SCR-space $M$ may be considered as a reflection
space in the sense of \Lit{LSSS}, i.e. if a `multiplication' on $M$ is
defined by $x\cdot y\colon= s_xy$ for all $x,y\in M$, the following rules
hold: $x\cdot x=x$, $\;x\cdot(x\cdot y)=y\;$ and $\;x\cdot(y\cdot
z)=(x\cdot y)\cdot(x\cdot z)$ for all $x,y,z\in M$. The SCR-spaces
form in various ways a category.
We prefer here the following notion (ignoring the Riemannian metrics on $M$ and $N$):

\Definition{} A CR-map $\phi\colon M\to N$ is called an SCR-map,
if $\phi(x\cdot y)=\phi\Kl x \cdot\phi(y)$ for all $x,y\in M$.
\Formend

If $\phi$ in addition is contractive we also call it a {\sl metric SCR-map}. In this
sense it is clear what (metric, isometric) SCR-isomorphisms,
SCR-automorphisms are. For instance, $G_M$ consists of
isometric SCR-automorphisms of $M$. Also, the universal covering
$\pi\colon\widetilde M\to M$ of an SCR-space $M$ has a unique structure
of an SCR-space such that locally $\pi$ is an isometric CR-diffeomorphism.

The following statement can be used to construct SCR-spaces.
\Lemma{SIGS} Let $M$ be a connected Hermitian CR-space with a base point $o$ and let
$H\subset I_M$ be a subgroup acting transitively on $M$.
Suppose that $\sigma\colon M\to M$ is a
diffeomorphism with $\sigma\Kl o=o$ and $\sigma^2=\id$
such that the following conditions are satisfied:
\0 $\sigma\circ H=H\circ\sigma$,
\1 the differential $d_o\sigma\in\5L(T_oM)$ is a linear isometry
with $(H_oM\oplus T_o^{rr}M)\;\subset\;\Fix(-d_o\sigma)$.\par
\noindent
Then $\sigma$ is a symmetry of $M$.
\Proof Fix $a\in M$ and choose $g,h\in H$
with $g\Kl a=o$ and $\sigma=h\,\circ\,\sigma\,\circ\,g$.
The claim follows from the identity $d_a\sigma=d_oh\,\circ\,d_o\sigma\,\circ\,d_a g$.
\qed

\Example{EXAP} Let $E,F$ be complex Hilbert spaces of finite
dimension. Suppose that $w\mapsto w^*$
is a conjugation of $F$ (i.e. a conjugate linear, involutive
isometry of $E$) and let
$\Phi\colon E\times E\to F$ be an Hermitian mapping with respect to
the conjugation $^*$ (i.e. $\Phi$ is sesqui-linear
and $\Phi(v,u)=\Phi(u,v)^*$ for all $u,v\in E$).
Set $$V\colon=\{w\in F:w+w^*=0\}\steil{and}
M\colon=\{(z,w)\in E\oplus F:w+w^*=\Phi(z,z)\}\,.$$
Then $M$ is a CR-submanifold of $E\oplus F$ with
$$\eqalign{T_aM=&\{(z,w)\in E\oplus F:w+w^*=\Phi(e,z)+\Phi(z,e)\}\cr
  H_aM=&\{(z,w)\in E\oplus F:w=\Phi(e,z)\}\cr}$$
for every $a=(e,c)\in M$. The group
$$\Lambda\colon=\big\{(z,w)\mapsto(z+e,w+\Phi(e,z)+\Phi(e,e)/2
+v):e\in E,v\in V\big\}$$
acts transitively and freely on $M$ by affine CR-diffeomorphisms. Therefore $M$ has a
group structure, a generalization of the Heisenberg group.
For $o\colon=(0,0)\in M$,
there exists a unique $\Lambda$-invariant Riemannian metric
on $M$ such that $T_oM=E\oplus V$ is the orthogonal sum of $E$
and $V$. By Lemma~\ruf{SIGS}, $M$ is a symmetric CR-manifold
-- for every $a=(e,c)\in M$, the corresponding symmetry $s_a$ is given by
$(z,w)\mapsto\big(2e-z,w+\Phi(2e,e-z)\big)$ and
$G_M=\Lambda\cup s_o\Lambda$, $G^0_M=\Lambda$.
The full group $\AUT(M)$ does not act properly on $M$
since it contains all transformations of the form
$(z,w)\mapsto(tz,t\overline t\[1w)$, $t\in\CC^*$.
\hfill\break
For every $\xi\in E$, the vector field $X$ on $E\oplus F$ defined
by $X_{(z,w)}=(\xi,\Phi(z,\xi))$ satisfies $X_a\in H_aM$ for all $a\in M$.
From this it is easily derived that the Levi form $L_o$ at $o\in M$ as
defined in \Ruf{LEVI} coincides with the Hermitian map $\Phi\colon E\times E\to F$
after the identification $H_oM=E$ and $(T_oM/H_oM)\otimes_\RR\CC\cong F$.
\qex
The next statement will be used later.
\Lemma{SPAT} Let $E,V\subset F,M,\Phi,\Lambda$ be as in Example \ruf{EXAP}. Assume that
$\Phi$ is non-degenerate in the following sense: For every $e\in E$ with
$e\ne0$ there exists $c\in E$ with $\Phi(e,c)\ne0$. Then
$$\Aff(M)=\Lambda\rtimes\GL(M)\;\;\4{and}$$
$$\GL(M)=\{(\eta\times\epsilon)\in\GL(E)\times\GL(V):
\Phi(\eta z,\eta z)=\epsilon\[1\Phi(z,z)\steil{for all}z\in E\}\,.$$
Furthermore, the group $I_M$ of all isometric CR-diffeomorphisms of $M$ is given by
$$I_M=\Lambda\rtimes\,\Gamma\subset\Aff(M)\,,\steil{where}
\Gamma\colon=\{(\eta\times\epsilon)\in\GL(M):\eta\;\4{unitary},\,
\epsilon\;\4{orthogonal}\}\,.$$
\Proof Let us start with an arbitrary real-analytic CR-diffeomorphism
$\phi$ of $M$ satisfying $\phi(o)=o$. Then $\phi$ extends to a holomorphicühic
map $\phi\colon U\to E\oplus F$ for a suitable open connected 
neighbourhood $U$ of $M$ in $E\oplus F$ (see e.g. \Lit{BER}, \S1.7). The differential 
$g\colon=d_o\phi\in\GL(E\oplus F)$ leaves $H_oM=E$
invariant and hence can be written
as operator matrix $g={\eta\,\alpha\choose0\,\epsilon}\in\GL(E\oplus F)$ with 
a linear operator $\alpha\colon F\to E$. Since also $T_oM=E\oplus V$ is
invariant under $g$, the operator $\epsilon\in\GL(F)$ must leave invariant
the subspace $V\subset F$, i.e. $\epsilon\in\GL(V)\subset\GL(F)$ and, in particular,
$\epsilon(w^*)=(\epsilon w)^*$ for all $w\in F$.
There exist holomorphic functions
$h\colon U\to E$ , $f\colon U\to F$ vanishing of order $\ge2$ at $o$
such that
$$\phi(z,w)\;=\;\big(\eta z+\alpha w+h(z,w),\epsilon w+f(z,w)\big)$$
for all $(z,w)\in M$. For every $z\in E$, $v\in V$ and 
$$w=w(z,v)\colon=v+\Phi(z,z)/2$$ we have $(z,w)\in M$ and hence
$$\epsilon\,\Phi(z,z)+f(z,w)+f(z,w)^*
=\Phi(\eta z+\alpha w+h(z,w),\eta z+\alpha w+h(z,w))\leqno{(*)}$$
for all $z\in E$ and $v\in V$. Comparing terms in $(*)$ we derive
$\epsilon\,\Phi(z,z)=\Phi(\eta z,\eta z)$ for all $z\in E$.
Now suppose that $\phi$ is affine, i.e. $f=0$ and $h=0$. Comparing 
terms in $(*)$ again we get
$\Phi(\alpha v,\eta z)=0$ for all $z\in E$ and $v\in V$. But then
the non-degeneracy of $\Phi$ implies $\alpha v=0$ for all $v\in V$,
i.e. $\alpha=0$. This proves that the groups
$\Aff(M)$ and $\GL(M)$ have the claimed forms.\hfill\break
Now suppose that $\phi\in I_M$ is an isometry. Since there is a unique
real-analytic structure on $M$ such that the Lie group $I_M$ acts as a 
real-analytic transformation group. Since the same holds
for the Lie subgroup $\Lambda$, these two structures must coincide, i.e. $\phi$ is 
real-analytic. Furthermore, $g=d_o\phi$ is a linear isometry of $E\oplus V$, i.e. 
$g={\eta\,0\choose0\,\epsilon}$ with $\eta$ unitary and $\epsilon$ orthogonal.
Together with $\epsilon\,\Phi(z,z)=\Phi(\eta z,\eta z)$ for all $z\in E$
this implies $g\in I_M$ and hence $\phi=g\in\GL(M)$
as a consequence of the Uniqueness Theorem \ruf{UUNI}.
\qed

For all CR-manifolds $M$ in Example \ruf{EXAP} with $\Phi\ne0$ the group
$G^0$ is nilpotent of nilpotency class $2=\kappa(M)$. For examples with nilpotent
groups of higher class compare section \ruf{sechs}. An explicit example
of class $3$ with lowest possible dimension is the following.
\Example{THRE} Set
$$M\colon=\big\{(z,w,v)\in\CC^3:\Im\Kl w
=z\overline z,\;\Im\Kl v=\Im(w\overline z)\big\}\,.$$
Then $M$ is a CR-submanifold of $\CC^3$ with CR-dimension 1 and CR-codimension
2 and every $(a,b,c)\in M$ induces an affine CR-automorphism of $M$ by
$$\big(z,w,v\big)\;\longmapsto\;\big(z+a,\,w+2i\overline az+b,\,
v+(2i\overline a^2-\overline b)z+(a+2\overline a)w+c\big)\,.\leqno{(*)}$$
Indeed, if we denote the right hand side of $(*)$ by ${\bf(z,w,v)}$ then
$$\eqalign{\Im({\bf w\overline z})\;&
=\;\Im(w\overline z+2iz\overline z\overline a+\overline zb
+2i\overline a^2z+\overline aw+\overline ab)\cr
&=\;\Im(v+(w-\overline w)\overline a-\overline bz+2i\overline a^2z+
\overline aw+c)\;=\;\Im({\bf v})\,.\cr}$$
An elementary calculation shows that
the transformations $(*)$ form a nilpotent Lie group $G^0$ acting
freely and transitively on $M$. In particular, $M$ has the structure of a
group with the product
$(a,b,c){\scriptstyle\,\odot\,}(z,w,v)\colon={\bf(z,w,v)}$
and the unit $o\colon=(0,0,0)$. There is a unique $G^0$-invariant
Riemannian metric on $M$ whose restriction
to the tangent space $T_oM=\CC\oplus\RR^2\subset\CC^3$
is the one inherited from $\CC^3$.
The transformation $s_o:(z,w,v)\mapsto(-z,w,-v)$
is a symmetry at $o$ by Lemma~\ruf{SIGS}.
Hence $M$ is a symmetric CR-manifold. It can be verified that $I_M=G_M=G^0\cup s_oG^0$
is the group of all CR-isometries. The action of $\AUT(M)$ is not proper
since this group contains all transformations of the form
$(z,w,v)\mapsto(tz,t^2w,t^3v)$, $t\in\RR^*$.\hfill\break
The group $Z$ of all translations $(z,w,v)\mapsto(z,w,v+c)$,
$c\in\RR$, is in the center
of $G^0$ and $M/Z$ is CR-isomorphic to the classical Heisenberg group
$\{(z,w)\in\CC^2:\Im\Kl w=z\overline z\}$ that already occurred in
Example \ruf{EXAP} in a slightly different form.
\qex

\KAP{vier}{Some examples derived from the sphere}

We start with the inclusion of the complex manifolds
$$\BB_n\subset\CC^n\subset\PP_n\,,\Leqno{SHSP}$$
where $\BB_n\colon=\BB\colon=\{z\in\CC^n:(z|z)<1\}$ is the euclidean ball
and $\PP_n=\PP_n(\CC)$ is the complex projective space
with homogeneous coordinates $[z_0,z_1,\dots,z_n]$. We identify
every $z\in\CC^n$ with the point $[1,z]\in\PP_n$. It is known
that the group $\Aut(\PP_n)$ of all biholomorphic automorphisms
of $\PP_n$ coincides with the group of all projective linear transformations $\PSL(n+1,\CC)$.
Furthermore,
$$\Aut(\BB)\;=\;\{g\in\Aut(\PP_n):g(\BB)=\BB\}\;=\;\PSU(1,n)\,.$$
The group $\Aut(\BB)\subset\Aut(\PP_n)$ has three orbits in $\PP_n$ --
the ball $\BB$, the unit sphere $S\colon=\partial\BB$ and the outer domain
$\DD\colon=\PP_n\backslash\overline\BB$. Actually, it is known
that again $$\Aut(\DD)\;=\;\{g\in\Aut(\PP_n):g(\DD)=\DD\}\;=\;\Aut(\BB)$$
holds. The spaces $\BB_n$, $\CC^n$ and $\PP_n$ are symmetric Hermitian manifolds
with constant holomorphic sectional curvature $<0$, $=0$ and $>0$ respectively.
Also, $\BB_n$ and $\PP_n$ are dual to each other in the sense of
symmetric Hermitian manifolds.

\Example{EXAM} The unit sphere
$S=\partial\BB=\{z\in\CC^n:(z|z)=1\}$
is a CR-submanifold of $\CC^n$, whose holomorphic tangent space at $a\in S$
is $H_a=\{v\in\CC^n:(a|v)=0\}$.
To avoid the totally real case $n=1$ let
us assume for the rest of the section that always $n>1$ holds.
Then $S$ is a minimal CR-manifold. It is known (compare f.i. \Lit{TH})
that for every pair $U,V$ of domains in $S$ and every CR-diffeomorphism
$\phi\colon U\to V$ there exists a biholomorphic transformation
$g\in\Aut(\PP_n)$ with $\phi=g|U$. In particular this implies
$$\AUT(S)\;=\;\{g\in\Aut(\PP_n):g(S)=S\}\;=\;\Aut(\BB)$$
and the maximal compact subgroups of $\AUT(S)$ are in one-to-one
correspondence to the points of $\BB$.\hfill\break
By restricting the flat Hermitian metric
of $\CC^n$, the  sphere $S$ becomes an Hermitian CR-manifold.
The unitary group $\U(n)$ coincides with $\{g\in\AUT(S):g(0)=0\}$,
acts transitively on $S$ by isometric
CR-diffeomorphisms and it is easy to see that actually
$I_M=\U(n)$ holds.
Moreover, $z\mapsto2(a|z)a-z$ defines a symmetry at $a\in S$
and $G_M=\{g\in\U(n):\det(g)=(\pm1)^{n-1}\}$.
The group $G^0=\SU(n)$ is simple whereas $I_M=\U(n)$ has
center $Z\cong\TT$. For every closed
subgroup $A\subset Z$ also $S/A$ is an SCR-space
with CR-dimension $n-1$ in a natural way -
for instance for $A=\{\pm\id\}$ we get
the real projective space $\PP_{2n-1}(\RR)$ and for $A=Z$
the complex projective space $\PP_{n-1}(\CC)$.~
\qex

The space of (projective) hyperplanes $L\subset\PP_n$ is again
a complex projective space of dimension $n$ on which the
group $\AUT(S)\subset\Aut(\PP_n)$ acts with three
orbits. These consist of all $L$ meeting $S$ in no, in precisely one
and in more than one point respectively.
Assume $L\cap S\ne\emptyset$ in the following and consider the domain
$W\colon=S\backslash L$ in $S$. Let $o\in W$ be the point with
the maximal distance from the hyperplane $L\cap\CC^n$.
We will see that $W$ has the structure of a symmetric
CR-manifold.
We claim that there exists a CR-isomorphic model $Q\subset\CC^n$
of $W$ in such a way that $Q$ is closed in $\CC^n$ and such
that every CR-diffeomorphism of $Q$ is the restriction of a complex
affine transformation of $\CC^n$, that is $\AUT(Q)=\Aff(Q)$. 
Indeed, choose a transformation
$g\in\Aut(\PP_n)$ with $g(L)\cap\CC^n=\emptyset$ and put $Q\colon=g(W)$.
We call $Q$ an {\sl affine model} of $W$.
Let us consider the two cases $L\cap\BB=\emptyset$
and $L\cap\BB\ne\emptyset$ separately.
\Example{EXAA} Let $U\colon=S\backslash L$ for a hyperplane $L$
with $L\cap\BB\ne\emptyset$, say $U=\{z\in S:z_1\ne0\}$ and
$o=(1,0,\dots,0)$. The group $\AUT(U)$ acts transitively on $U$
 and has compact isotropy subgroup $\U(n\]4-\]51)$ at $o$.
Therefore, $U$ is a symmetric CR-manifold with
$$I_U=\AUT(U)\;\cong\;\U(1,n\]4-\]51)$$
and the symmetry $\rho=s_o$ at $o$ is given by
$\rho(t,v)=(t,-v)$ for all $(t,v)\in\CC^{1+(n-1)}$.
An affine model is
$$R\colon=\{z\in\CC^n:(\rho z|z)=1\}\;=\;\{(t,v)\in
\CC^{1+(n-1)}:t\overline t-(v|v)=1\}\Leqno{DUAL}$$
with $(t,v)\mapsto (1/t,v/t)$ a CR-diffeomorphism $U\to R$.
The universal covering $\widetilde R$ of $R$ is again
a symmetric CR-manifold and can be realized via $(s,v)\mapsto(\exp\Kl s ,v)$ as
(see also Example \ruf{SPAT})
$$\widetilde R=\{(s,v)\in\CC^{1+(n-1)}:\exp(s+\overline s)-(v|v)=1\}\,.\eqno{\qex}$$
\Example{EXAX} Let $V\colon=S\backslash\{a\}$ for some point $a\in S$,
say $a\colon=(0,\dots,0,1)$ and
hence $o=-a$. Then $V$ is a cell in $S$ and
the Cayley transform $(v,t)\mapsto\big(\sqrt2\,(v/(1-t),(1+t)/(1-t)\big)$
for all $(v,t)\in\CC^{(n-1)+1}$, $t\ne1$, defines a
CR-diffeomorphism of $V$ onto the affine model
$$N\colon=\big\{(v,t)\in\CC^{(n-1)+1}:t+\overline t=(v|v)\big\}\,.$$
This SCR-space occurs already in Example \ruf{EXAP}
and as a consequence of \ruf{SPAT} we have
$$I_N\;=\;G_{\]3N}^0\[3\rtimes\,\U(n\]4-\]51)\steil{and}\AUT(N)\;=\;\Aff(N)\;
=\;G_{\]3N}^0\[3\rtimes\,\big(\RR^+\times\U(n\]4-\]51)\big)\,,$$
where the groups $\U(n\]4-\]51)$ and $\RR^+$ act on $N$ by
$(v,t)\mapsto(\epsilon v,t)$ and $(v,t)\mapsto(sv,s^2t)$ respectively.
Consider for every $s\in\RR^+$, the central subgroup
$\Gamma_s\colon=\{(v,t)\mapsto(v,t+ins):n\in\ZZ\}$ of $G_N$.
Then the quotient manifold $N_s\colon=N/\Gamma_s$ is an SCR-space diffeomorphic
to $\CC^n\times\TT$. But since the groups
$\Gamma_s$ and $\Gamma_{\tilde s}$ are not conjugate in
$\AUT(N)$ for $s\ne\tilde s$ the manifolds $N_s$ and $N_{\tilde s}$
are not isomorphic as CR-manifolds. We obtain a continuous family of symmetric CR-manifolds
that are pairwise non-isomorphic even in the category of CR-manifolds.
\qex
In analogy to \Ruf{SHSP} we have inclusions
$$U\subset V\subset S\,.\Leqno{INCL}$$
But in contrast to \Ruf{SHSP} $U$ is not `a bounded domain' (meaning
relatively compact) in the cell $V$.
In all three cases $M=U,V,S$, the center $Z$
of $I_M$ is either $\TT$ or $\RR$ and the quotient CR-manifold
$M/Z$ is $\BB_{n-1}$, $\CC^{n-1}$ and $\PP_{n-1}$ respectively.
We call the symmetric
CR-manifold $R\cong U$ the {\sl dual unit sphere} in $\CC^n$ (compare
also the discussion at the end of section \ruf{sechs}).
We remark that the action of $\TT$ on $U$ given by
scalar multiplication has as quotient the complex
manifold  $\CC^{n-1}$ which is not biholomorphically equivalent
to the bounded domain $\BB_{n-1}=U/Z$. 

The group $G_M$ is simple in case $M=U,S$ and is nilpotent
in case $M=V$. Also, the isotropy subgroup $K$ at a point $o\in M$
acts irreducibly on the tangent space $T_aM$ if $M=U,S$
and $K$ is a finite group in the third case. In all three cases
the group $I_M$ acts transitively on the subbundle
$\{v\in H_aM:a\in M,\,\|v\|=1\}$ of the tangent bundle $TM$.
In particular, $M$ is an SCR-space of constant holomorphic sectional curvature.

\KAP{fuenf}{A canonical fibration}

Let again $M$ be an SCR-space and let $o\in M$ be a
fixed point, called base point in the following.
Let $G=G_M$ and $K=\{g\in G:g\Kl o=o\}$ be as before.
Then $s_o$ is in the center of $K$ and $\sigma(g)\colon= s_og s_o$
defines an involutive group automorphism $\sigma$ of $G$. Therefore
$K$ is contained in the closed subgroup $\Fix(\sigma)\subset G$.
Let $L$ be the smallest {\sl open} subgroup of $\Fix(\sigma)$ containing $K$.
Then $s_o$ is contained in the center of $L$ and $s_a= s_o$ for
all $a\in L\Kl o $. Identify as before
$M$ with the homogeneous space $G/K$ and put $N\colon=G/L$.
Then we have canonical fibre bundles
$$G\buildrel \mu\over\longrightarrow M\buildrel\nu\over\longrightarrow N$$
defined by $g\mapsto gK\mapsto gL$. The typical fibres are
$K$ for $\mu$ and the connected homogeneous space  $L/K$ for $\nu$.
The following statement follows directly from the definition of $\sigma$.
\Lemma{LIFT} The fibration $\mu$ satisfies $\mu\circ\sigma= s_o\circ \mu$.
\Formend

Because of Lemma~\ruf{LIFT}, $\sigma$ can be seen as a lifting of $s_o$ via $\mu$.
On the other hand, $s_o$ can be pushed forward via $\nu$:

\Proposition{KLAS}  For every $c\in N$, there exists a unique involutive
diffeomorphism $\, s_c\colon N\to N$ such that
$\nu\circ s_a= s_c\circ\nu$ for all $a\in M$ with $\nu\Kl a=c$.
The differential $d_a\nu$ has kernel $T_a^+M$ in $T_aM$ and hence
induces a linear isomorphism from $T_a^-M$ onto $T_cN$.
Every $s_c$ has $c$ as an isolated fixed point.
Furthermore, $N$ is simply connected if $M$ is simply connected.
\Proof Since the fibration $\nu$ is $G$-equivariant we have to establish
the map $s_c$ only for the base point $c:=\nu\Kl o $ of $N$.
But then $s_c$ is given by $gL\mapsto\sigma(g)L$ since $s_o$ can be
identified with the map $gK\mapsto\sigma(g)K$ of $M$. The tangent space
$T_oF$ of the fiber $F\colon=\nu^{-1}(c)=L\Kl o $ at $o$ is $T_o^+M$ and
hence is the kernel of the differential $d_o\nu$. Consequently,
the differential $d_c s_c$ is the negative identity on $T_cN$ and
hence $c$ is an isolated fixed point of $s_c$.
Now suppose that $M$ is simply connected
and denote by $\alpha\colon H\to G^0$ the universal covering group
of $G^0$. Then the subgroup $\alpha^{-1}(K\cap G^0)$ of $H$ is connected
and hence by the construction of $L$ also the subgroup
$\alpha^{-1}(L\cap G^0)$ is connected, i.e. $N=G^0/(L\cap G^0)$ is simply connected.
\qed

The manifold $N$ in Proposition \ruf{KLAS} together with
all involutive diffeomorphisms $s_c$, $c\in N$, is a
symmetric space in the sense of \Lit{LOSS}.
An interesting case occurs when $N$ has the structure of a symmetric
CR-manifold in such a way that
{\parskip0pt
\0 $\nu$ is a CR-map and $M$, $N$ have the same CR-dimension.
\1 $\nu$ is a partial isometry, i.e. the restriction of $d_a\nu$ to
  $T_a^-M$ is an isometry for every $a\in M$.
\1 $s_c$ is a symmetry at $c$ for every $c\in N$.\par}
\noindent
We say that the SCR-space $M$ has {\sl symmetric reduction} $N$ if the
properties (i)  -- (iii) are satisfied. For every $g\in L$ and
$a=g\Kl o $ we have the commutative diagram

$$\diagram{T_o^-M&\mapright{\displaystyle d_o\Phi_g} &T_a^-M\cr
\mapdown{\displaystyle\hskip-23pt d_o\nu}&&\mapdown{\[4\displaystyle d_a\nu}\cr
T_cN&\mapright{\displaystyle d_c\]1\Psi_g} &T_cN\cr}$$
where for better distinction we denote by $\Phi_g$ the diffeomorphism
of $M$ given by $g$ and by $\Psi_g$ the corresponding diffeomorphism of $N$ (that
is, $\Phi_g(hK)=ghK$ and $\Psi_g(hL)=ghL\;)$.
Put $H_cN\colon=d_o\nu(H_oM)$ and give it the complex structure for which
$d_o\nu\colon H_oM\to H_cN$ is a complex linear ismorphism. Furthermore,
endow $T_cN$ with the Riemannian metric
for which $d_o\nu\colon T_o^-M\to T_cN$ is an isometry.
Then the existence of a $G$-invariant almost CR-structure on $N$
with property (i)  is equivalent
to the condition that all operators $d_c\]1\Psi_g$, $g\in L$,
leave the subspace $H_cN$ invariant and are complex linear there.
In the same way, a $G$-invariant Riemannian metric on $N$ with property (ii)
exists if and only every $d_c\]1\Psi_g$, $g\in L$, is an isometry of $T_cN$.
This happens for instance (after possibly changing the metric of $M$)
if the group $L$ is compact,
or more generally, if the linear group $\{d_c\]1\Psi_g:g\in L\}$ is
compact.

We notice that as a consequence of Lemma
\ruf{SIGS}, condition (iii) is automatically satisfied if (i), (ii) hold. Also,
in case that for the fibration $\nu\colon M\to N$ there exists a
Riemannian metric on $N$ with the property (ii) and such that in addition
every $\nu$-fibre is a symmetric Riemannian manifold, the space
$M$ is a bisymmetric space in the sense of \Lit{KANT}.

The following sufficient condition for the existence of a symmetric
reduction is easily seen, we leave the proof to the reader.
\Lemma{AKA} Suppose that the subgroup $\{g\in I_M:\nu\circ g=\nu\}$ acts
transitively on some $\nu$-fibre. Then this group acts transitively
on every $\nu$-fibre and $M$ has a symmetric reduction.
\qex

\medskip
Not every symmetric CR-space has a symmetric reduction. Consider
for instance Example \ruf{THRE}. Then $L^0$ is the subgroup
of all $(0,b,0)\in M$ with $b\in\RR$, and $N\colon=G/L$ can be identified
with $\CC\times\RR$ in such a way that 
$\nu(z,w,v)=\big(z,\Re(v)-3\Re(w)\Re(z)\big)$.
Furthermore, the action of $L^0$ on $\CC\times\RR$ is given by
$(z,t)\mapsto\big(z,t-4\[1b\[1\Re(z)\big)$, $b\in\RR$.
This implies that there cannot exist any $G$-invariant Riemannian metric
on $N$. Also, there is no CR-structure on $N$ satisfying property (i).

\KAP{fuenfa}{A construction principle}

In this section we give a Lie theoretical construction of symmetric CR-spaces
such that every SCR-space can be obtained is this way.
We start with an arbitrary connected Lie group $G^0$ together with
an involutive group automorphism $\sigma$ of $G^0$.
Then there is a Lie group $G$ with connected identity component
$G^0$ and an element $s\in G$ with $G=G^0\,\cup\,s\[2G^0$
and $\sigma(g)=sgs$ for all $g$.
Let $\7g$ be the Lie algebra of $G$ and denote by $\tau\colon=\Ad(s)$
the Lie algebra automorphism of $\7g$ induced by $\sigma$ (here and in the
following $\Ad$ always refers to the group $G$). Put
$$\7l\colon=\Fix(\tau)\steil{and}\7m\colon=\Fix(-\tau)\Leqno{FIXT}$$
Then $\7l$ is a Lie subalgebra of $\7g$ and $\7m$ is
a Lie triple system, see \Lit{LOSS}. For every $g\in\Fix(\sigma)$,
the decomposition $\7g=\7l\oplus\7m$
is invariant under $\Ad(g)$.
\hfill\break
Now choose a compact subgroup $K\subset\Fix(\sigma)$, an $\Ad(K)$-invariant
Riemannian metric on $\7g$ and a linear
subspace $\7h\subset\7m$ together with a complex structure $J$ on $\7h$
satisfying the following properties:
{\parskip0pt
\0 $\|Jx\|=\|x\|$ for all $x\in\7h$.
\1 $K$ contains the element $s$.
\1 $\Ad(g)$ leaves the subspace $\7h$ invariant
and commutes there with $J$ for every $g\in K$.\par}
\noindent
Notice that $K=\{s,e\}$ with arbitrary $\7h\subset\7m$ and
arbitrary $J$ always is an admissible choice.
Also, if the compact group $K$ has been chosen,
every closed subgroup of $K$ containing
$s$ is again an admissible choice.

Since $K$ is compact there exists an
$\Ad(K)$-invariant decomposition $\7l=\7k\oplus\7n$,
where $\7k$ is the Lie algebra of $K$. With $\7p\colon=\7n\oplus\7m$ therefore
we get the $\Ad(K)$-invariant decomposition $\7g=\7k\oplus\7p$.
Consider the connected homogeneous $G$-manifold $M\colon=G/K$ and declare
$o\colon=K\in M$ as base point. In the following we identify the tangent space $T_oM$
in the canonical way with the Hilbert space $\7p$. Denote by
$\Phi_g$ the diffeomorphism of $M$ induced by $g$ , that is
$$\Phi_g\colon M\to M,\quad hK\longmapsto ghK\,,$$
(we do not require here
that the $G$-action is effective, this could be easily achieved by reducing
out the kernel of ineffectivity from the beginning). Then for every $g\in K$,
the differential $d_o\Phi_g$ is nothing but the restriction of
$\Ad(g)$ to $\7p$ and hence there exists a unique $G$-invariant almost CR-structure
with $H_oM=\7h$ and also a unique $G$-invariant Riemannian metric on $M$
extending the given Hilbert norm of $T_oM=\7p$. In particular, $s_o\colon=\Phi_s$
is an involutive isometric diffeomorphism of $M$ with fixed point $o$.
Clearly, $s_o$ is a symmetry of $M$ at $o$ if and only if $H_o^{-1}M\subset\7m$,
where $H_o^{-1}M\subset T_oM$ is the subspace defined in section \ruf{zwei}.
A more convenient condition for this is given by the following statement.
\Proposition{} Let $\7a$ and $\7b$ be the Lie subalgebras of $\7g$ generated
by $\7m$ and $\7h$ respectively. Then $\7a=[\7m,\7m]\oplus\7m$ holds,
$\7a$ is an $\Ad(K)$-invariant ideal of $\7g$ and
\0 $M$ is a minimal symmetric CR-manifold with symmetry $s_o$ at $o$
if and only if $\7g=\7k+\7b$.
\1 In case $s_o\colon=\Phi_s$ is a symmetry of $M$ at $o$, the weaker 
condition $\7g=\7k+\7a$ holds.
\Proof First notice that $[\7m,\7m]\subset\7l$, $[\7l,\7m]\subset\7m$ by \Ruf{FIXT}
and hence that $\7a=[\7m,\7m]\oplus\7m$ holds. Obviously, $\7a$ is invariant
under $\ad(\7l)$ as well as $\ad(\7m)$, i.e. $\7a$ is an ideal in $\7g$.
Now suppose that $s_o$ is a symmetry of $M$ at $o$. Then $M$ is a symmetric CR-manifold
by the transitivity of the group $G$ and $\7g=\7k+\7a$ follows as in the proof
of Proposition \ruf{SYMM}. Now suppose that
in addition that $M$ is minimal as an almost CR-manifold.
Define inductively $\7h^k\colon=\7h^{k-1}+[\7h,\7h^{k-1}]$ and $\7h^0=0$.
Then $\7h^k/\7h^{k-1}$ is isomorphic to $H_o^kM$ and $\7g=\7k+\7h^k$ for $k$ sufficiently large,
i.e. $\7g=\7k+\7b$. Conversely, suppose that $\7g=\7k+\7b$ holds.
Then $M$ is minimal and the differential of $s_0$ is the negative identity
on $\7h=H_oM$, i.e. $s_o$ is a symmetry at $o$.\qed

As an illustration of the construction principle fix integers
$p,q\ge0$ with $n\colon=p+q\ge2$ and set $G\colon=\SU(n)$.
Then the corresponding Lie algebra $\7{g=su}(n)$ is a real Hilbert
subspace of $\CC^{n\times n}$. Write every $g\in\CC^{n\times n}$ in the
form ${a\[1b\choose c\[1d}$ with $a,b,c,d$ matrices of sizes
$p\times q, p\times q, q\times p$, $q\times q$ respectively and
denote by $\sigma$ the automorphism of $G$ defined by ${a\[1b\choose c\[1d} 
\mapsto{\;a\;-b\choose-c\;\,d}$. Fix a closed subgroup 
$K\subset F\colon=\Fix(\sigma)$
and put $M\colon=G/K$ with base point $o\colon=K\in M$. Identify
$T_oM$ with the orthogonal complement $\7p$ of $\7k$ in $\7g$
and put $H_oM\colon=\7m\colon=\{{0\[1b\choose c\[20}\in\7p\}\approx\CC^{p\times q}$
with complex structure defined by ${0\[3b\choose c\[40}\mapsto
{\,\[90\[{15}ib\choose-ic\[60}$. These data give a unique $G$-invariant
Hermitian metric and a unique $G$-invariant almost CR-structure on $M$
with $gK\mapsto\sigma(g)K$ a symmetry at $o\in M$. It is easily seen
that $\7m+[\7m,\7m]=\7g$ as well as the integrability condition 
\Ruf{INTE} hold, i.e. $M$ is
a compact minimal symmetric CR-manifold with symmetric reduction
$N\colon=G/L$. Here $N$ is the 
Grassmannian $\GG_{p,q}$ of all linear subspaces of dimension $p$ in
$\CC^n$ and in particular is a symmetric Hermitian space.
If we replace $G=\SU(n)$ by the group $G^d\colon=\SU(p,q)$ and define
$\sigma$ by the same formula, then $L=\Fix(\sigma)$ remains unchanged and for
every compact subgroup $K\subset L$ we get the two minimal symmetric
CR-manifolds $M=G/K$ and $M^d\colon=G^d/K$ which we call dual to
each other. In particular, $N^d\colon=G^d/L$ is a bounded symmetric
domain and is the dual of the Grassmannian 
$\GG_{p,q}$ in the sense of symmetric Hermitian spaces.

\KAP{sechs}{Integrability and complexifications}
Assume that $M$ is an arbitrary CR-space with base point
$o\in M$, not necessarily symmetric to begin with.
Assume that $G$ is a Lie group acting smoothly and transitively on $M$
by CR-diffeomorphisms. Let $K\colon=\{g\in G:g\Kl o=o\}$ be the
isotropy subgroup at $o$ and denote by $\7k\subset\7g$ the
corresponding Lie algebras. Then the canonical map $\theta\colon\7g\to T_oM$
is surjective and has $\7k$ as kernel. Choose a linear subspace
$\7h\subset\7g$ such that $\theta\colon\7h\to H_oM$ is a linear
isomorphism. Then there is a unique complex structure $J$ on
$\7h$ making $\theta|_{\7h}$ complex linear.
It is clear that $\7k\oplus\7h=\theta^{-1}(H_oM)$ does
not depend on the choice of $\7h$.

Let $\8g\colon=\7g\oplus i\7g$ be the complexification
of $\7g$ and denote for linear subspaces of $\7g$ its complex linear span
by the corresponding {\sl boldface letter,}
that is for instance $\8a=\7a\oplus i\7a$ in case of $\7a$.
The complex structure $J$ of $\7h$ extends in a unique way to a
complex linear endomorphism $\3J$ of $\8h$. Denote by
$\8h^\pm$ the eigenspaces of $\3J$ in $\8h$ to the
eigenvalues $\pm i$: $\8h^\pm=\{Jx\pm ix:x\in\7h\}$.
Put $\8l\colon=\8k\oplus\8h^-$ for the following.
This space does not depend on the choice of $\7h$.
The composition of the canonical maps $\7h\hookrightarrow\8h\to\8h/\8h^-$
induces a complex linear isomorphism $\7h\,\cong\,\8h/\8h^-\cong\,\8h^+$.
As a consequence, $\7g/\7k\hookrightarrow\8g/\8l$
realizes $T_oM=\7g/\7k$ as a linear subspace of the complex vector
space $\8g/\8l$ in such a way that there $H_oM=T_oM\cap iT_oM$ holds.
This property will be the key in the proof of Proposition \ruf{INBO}.

\Proposition{INBE} The CR-structure of $M$ is integrable if and only if
$\8l$ is a Lie algebra.
\Proof Let $\3TM$ be the complexified tangent bundle of $M$ and
denote by $\8V=\7V\oplus i\7V$ the complex Lie algebra of all smooth
complexified vector fields on $M$, i.e. of all smooth sections
$M\to\3TM$. Then $\8K\colon=\{X\in\8V:X_o=0\}$ is
a complex Lie subalgebra of $\8V$.
The almost CR-structure of $M$ gives a complex subbundle
$\3H^{0,1}M\subset\3TM$. Denote by $\8H^{0,1}\subset\8V$
the linear subspace of all vector fields
of type $(0,1)$, i.e. $X_a\in\3H^{0,1}_a\]2M$ for all $a\in M$.
The integrability conditon for $M$ is equivalent
to $\8H^{0,1}$ being a Lie subalgebra of $\8V$.
The action of $G$ on $M$ induces a Lie homomorphism
$\7g\to\7V$ that uniquely extends to a complex linear homomorphism
$\phi\colon\8g\to\8V$. The assumption that $G$ acts by CR-diffeomorphisms
implies $\big[\phi(\8g),\8H^{0,1}\big]\subset\8H^{0,1}$. For 
$\8L\colon=\8K+\8H^{0,1}$ we have $\8l=\phi^{-1}(\8L)$.

\noindent
Now suppose that $M$ is integrable. We claim that $\8l$ is a Lie algebra
and consider arbitrary vector fields $X,Y\in\phi(\8l)$. It is enough to
show for $Z\colon=[X,Y]$ that $Z_o\in\3H_o^{0,1}\]2M$ holds, 
i.e. that $Z$ is contained in $\8L$. Write $X=X'+X''$,
$Y=Y'+Y''$ with $X',Y'\in\8K\,$ and $X'',Y''\in\8H^{0,1}$. Then we have
$$[X',Y'']=[X-X'',Y'']\equiv[X,Y'']\steil{and hence}$$
$$Z\;\equiv\;[X',Y'']+[X'',Y']\;\equiv\;[X,Y'']+
[X'',Y]\;\equiv\;0\steil{mod}\8L\;.\leqno{(*)}$$
Conversely, suppose that $\8l$ is a Lie algebra.
Since $G$ acts transitively on $M$ we have $\8L=\8K+\phi(\8l)$.
We claim that $M$ is integrable.
Consider arbitrary vector fields $X,Y\in\8H^{0,1}$ and write $X=X'+X''$,
$Y=Y'+Y''$ with $X',Y'\in\8K\,$ and $X'',Y''\in\phi(\8l)$. We have to show that
$Z\colon=[X,Y]$ is contained in $\8H^{0,1}$. Since $G$
acts transitively on $M$ and leaves $\8H^{0,1}$ invariant
it is enough to show that $Z_o\in\3H_o^{0,1}\]2M$ holds, i.e. that $Z\in\8L$.
But this follows as in $(*)$.\qed
\Corollary{INBI} Suppose that $M$ in \ruf{INBE} is symmetric with $G\colon=G_M$.
Then, if $\7h\subset\7g$ is chosen to be $\Ad(K)$-invariant, $M$
is integrable if and only if $[\8h^-,\8h^-]\subset\8k$.
\Proof By the choice of $\7h$ we have $[\7k,\7h]\subset\7h$
and hence $[\8k,\8h^-]\subset\8h^-$.
The involution $\Ad(s_o)$ of $\7g$ extends to a complex linear
involution $\tau$ of $\8g$ with $\8h\subset\Fix(-\tau)$.
Therefore, $\8l$ is a Lie algebra if and only if 
the inclusion $[\8h^-,\8h^-]\subset\8l$ holds, that is
$[\8h^-,\8h^-]\;\subset\;\8l\cap\Fix(\tau)=\8k$.
\qed
We remark that \ruf{INBE} and \ruf{INBI} remain valid for
$\8h^+$ in place of $\8h^-$. 
\Proposition{INBO} Let $M=G/K$ be a homogeneous CR-manifold as in
Proposition \ruf{INBE}. Suppose, there exist complex Lie groups
$\3L\subset\3G$ with Lie algebras $\8l\subset\8g$, where $\8l$
and $\8g=\7g\oplus i\7g$ are as before.
Suppose furthermore that $G$ can be realized as real Lie subgroup
$G\subset\3G$ in such a way that the corresponding injection
$\7g\to\8g$ is the canonical one and such that $G\3L$ is locally
closed in $\3G$.
Then, if $\3L\cap G=K$ holds and if $\3L$ is closed in $\3G$, $gK\mapsto g\3L$
realizes $M$ as a locally closed generic CR-submanifold of the
homogeneous complex manifold $\3M\colon=\3G/\3L$.
\Proof The assumptions guarantee that $M$ is imbedded in $\3M$ as a locally
closed real-analytic submanifold with $H_oM=T_oM\cap iT_oM$ in $T_o\3M$.
The result follows since $M$ is a $G$-orbit in $\3M$.
\phantom{...}\qed
In general, the CR-submanifold $M$ is not closed in $\3M$.
For instance, if $M$ is a bounded symmetric domain and $G=\Aut(M)$
is the biholomorphic automorphism group of $M$ then $\3M$ can be chosen to be the 
corresponding compact dual symmetric Hermitian manifold which
contains $M$ as an open subset. For the sphere $M\colon=\partial\BB_n\subset\CC^n$
and $G=\AUT(M)$ we may chose $\3M=\PP_n$. On the other hand, for the same
sphere $M=\partial\BB_n$
but $G=\U(n)$ we may obtain for $\3M$ the domain
$\CC^n\backslash\{0\}$ in $\CC^n$ -- but also the Hopf manifold
$\Quot{\CC^n\backslash\{0\}}{\alpha^{\ZZ}}$ for some complex number $\alpha$
with $|\alpha|>1$.

In the following we illustrate the statements \ruf{INBE} -- \ruf{INBO}
by various examples. Denote by $\sigma$
the inner automorphism of $\CC^{n\times n}$ given by
$(a_{ij})\mapsto(\Kl{\]5-\]51}^{i+j}a_{ij})$.
\Example{VORH} On the contrary to symmetric Hermitian spaces,
the CR-structure of a symmetric CR-space does not need to be integrable. 
For $n\ge3$ let $M\subset\CC^{n\times n}$ be the
nilpotent subgroup of all unipotent lower triangular matrices, i.e. of
all $a=(a_{ij})$ with $a_{ii}=1$ and $a_{ij}=0$ if $i<j$. Then for
the identity $e\in M$, the tangent space $T_eM$ will be identified
with the nilpotent algebra $\7g$ of all strictly lower triangular matrices.
Denote by $G$ the group generated by all left multiplications with
elements from $M$ and denote the restriction of $\sigma$ to $M$ by the
same symbol. Then
$G=G^0\cup\sigma G^0$ acts transitively on $M$ and there
exists a unique $G$-invariant Riemannian metric on $N$ which
coincides on $T_eM$ with the one inherited from $\CC^{n\times n}$.
Also there is a unique $G$-invariant almost CR-structure on $M$
with $$H_eM=\7h\colon=\{a\in\7g:a_{ij}=0\;\;\4{if}\;j\ne i+1\}$$
and complex structure on $\8h$ inherited from $\CC^{n\times n}$.
With this structure $M$ is symmetric and minimal. Because
of $[\8h^-,\8h^-]\ne0$ the CR-structure is not integrable.
\qex

\Example{SYIN} Let $n>d\ge1$ be fixed integers with $d\le n/2$ and denote by $\7g$
the space of all matrices in $\CC^{n\times n}$ having the form \Ruf{TETT}.
Also, denote by $\7h\subset\7g$ the subspace of all matrices with%
\vadjust{\midinsert
$${\pmatrix{
0&\cr
z_1&0\hbox to 0pt{\vbox to 0pt{\hskip40pt\hbox{\gros 0}\vss}\hss}\cr
\alpha_1&\overline z_1&0\cr
z_2&-\alpha_1&z_1&0\cr
\alpha_2&\overline z_2&\alpha_1&\overline z_1&0\cr
z_3&-\alpha_2&z_2&-\alpha_1&z_1&0\cr
%\alpha_3&\overline z_3&\alpha_2&\overline z_2&\alpha_1&\overline z_1&0\cr
\vdots&\vdots&\vdots&\vdots&\vdots&&\;\ddots\cr
}}\Leqno{TETT}$$
\centerline{\vbox{\klein\hsize=5cm\noindent {\bf Example \ruf{SYIN}:} All
$z_k$ arbitrary complex numbers,
all $\alpha_k$ purely imaginary.}}
\endinsert}
$z_k=\alpha_j=0$ for all $k>d$ and all $j$.
A simple calculation shows that $\7g$ is a real Lie subalgebra
of $\CC^{n\times n}$ and that $\7h$
generates $\7g$ as Lie algebra. Identifying $z=(z_1,\dots,z_d)\in\CC^d$
in the obvious way with the corresponding matrix in $\7h$ we get
a complex structure $J$ on $\7h\subset\Fix(-\sigma)$.
For all $x,y\in\7h$ the identity $[Jx,y]+[x,Jy]=0$ is easily verified.
$M\colon\,=\,\exp(\7g)$ is a closed
nilpotent subgroup of $\GL(n,\CC)$ invariant under $\sigma$.
Precisely as in Example \ruf{VORH} $M$
becomes a symmetric minimal CR-manifold. But this time
$[\8h^-,\8h^-]=0$ holds, that is, $M$ is integrable.
Furthermore, $\kappa(M)=[(n-1)/(2d-1)]$
and $M$ has CR-dimension $d$.

Proposition \ruf{INBO} gives a prescription for a generic
embedding of $M$.
Let $\3G$ be the connected, simply connected complex Lie group
with Lie algebra $\8g=\7g\oplus i\7g$. Then $\exp\colon\8g\to\3G$ is biholomorphic
and in particular, $\3L\colon=\exp(\8h^-)$ is a closed
abelian complex subgroup of $\3G$ --
notice that we have $\8k=0$ in this case.
Now, $M$ embeds in the canonical way into the complex
manifold $\3M=\3G/\3L$ which is biholomorphic
to the complex vector space $\8g/\8h^-$. Easy for explicit calculations is the
case $n$ odd -- then $\7g\cap i\7g=0$ holds
in $\CC^{n\times n}$ and we can realize the complexification
$\8g$ within the complex Lie algebra $\CC^{n\times n}$.
The commutative subalgebra $\8h^-$ then consists of all matrices obtained
from \ruf{TETT} by keeping all $\overline z_1,\overline z_2,\dots,\overline z_d$
and replacing all other entries (including all $z_k$) by $0$.
For $n=3$, $d=1$ we find the realization
$$M\cong\{(z,w)\in\CC^2:w+\overline w=z\overline z\}$$
which is the classical Heisenberg group, compare Example \ruf{EXAP}. For $n=5$, $d=1$ one can show
$$\eqalign{M\cong\Big\{(z,w,v_1,v_2,u)\in\CC^5:w+\overline w=z\overline z,
\;\;v_1-\overline v_2&=z\overline z(z-\overline z)/6+\overline wz,\cr
u+\overline u&=w\overline w+(z\overline v_1+\overline zv_1)+z\overline zz\overline z/4
\Big\}\cr}$$
which is a symmetric CR-manifold with CR-dimension 1, CR-codimension 4 and $\kappa(M)=4$,
compare also Example \ruf{THRE}. The symmetry at the origin is given by
$(z,w,v_1,v_2,u)\mapsto(-z,w,-v_1,-v_2,u)$.
\qex\Formend

\Example{SPAT} Fix an integer $k>1$ and consider in $\CC^2$
the connected CR-submanifold
$$M\colon=\big\{(s,v)\in\CC^2:|s|^{2k}-|v|^2=1\big\}$$
which is a $k$-fold cover of the symmetric CR-manifold $R$
in Example \ruf{EXAA} via the map $(s,v)\mapsto(s^k,v)$.
Therefore also $M$ is a symmetric CR-manifold and $G^0_M$ is
a $k$-fold covering group of $G^0_R=SU(1,1)\cong\SL(2,\RR)$.
Denote by $\7g$ the Lie algebra of $G_M$. Then
it is known that for $\8g=\7g\oplus i\7g$ there does not exist
any complex Lie group $\3G$ into which $G$ admits an embedding
induced by the canonical injection $\7g\hookrightarrow\8g$.
Therefore the conclusion of Proposition \ruf{INBO} cannot hold
for this example.\hfill\break
The group $G^0_M$ consists of all transformations
$$(s,v)\mapsto\big((as^k+bv)^{1/k},\overline bs+\overline av\big)\,,$$
where $a,b\in\CC$ satisfy $a\overline a-b\overline b=1$. It follows that
the action of $G_M$ does not extend to all of $\CC^2$.
But it extends to the domain
$$D\colon=\big\{(s,v)\in\CC^2:|v|<|s|^k\big\}\;=\;\RR^+M$$
on which the group $\RR^+\!\!\times\!G^0_M$ acts transitively and
freely.
\qex

\KAP{sieben}{CR-manifolds derived from bounded symmetric domains}
Suppose that $E$ is a complex vector space of dimension $n$
and that $D\subset E$ is a bounded symmetric domain.
Then $\Aut(D)$ is a semi-simple Lie group and
at every point of $D$ the corresponding isotropy subgroup is a
maximal compact subgroup (see f.i. \Lit{HELG}).
It is known that there
exists a complex norm $\|\!\cdot\!\|_\infty$ on $E$
such that $D$ can biholomorphically be realized as the open unit ball
$$D=\{z\in E:\|z\|_\infty<1\}\Leqno{REAL}$$
with respect to this norm and that any two realizations of this
type are linearly equivalent. In this realization
the isotropy subgroup at the origin is linear, i.e.
$$\{g\in\Aut(D):g(0)=0\}\;=\;\GL(D)\,.$$
Moreover, $\GL(D)$ is compact. Therefore 
there exists a $\GL(D)$-invariant complex Hilbert
norm $\|\!\cdot\!\|$ on $E$ that we fix for the sequel
and hence consider $E$ as a complex Hilbert space in the following.
We will also always assume
that $D$ is given in the form \Ruf{REAL}.
For shorter notation we use for the whole section the abbreviation
$$\Gamma\colon=\Aut(D)^0\steil{and}K\colon=\GL(D)^0\,.$$
As a generalization of \ruf{SHSP} there exists a compact complex manifold
$Q$ with
$$D\subset E\subset Q\steil{and}\Aut(D)=\{g\in\Aut(Q):g(D)=D\}\,.\Leqno{COMP}$$
$Q$ is the dual of $D$ in
the sense of symmetric Hermitian manifolds and $\Aut(Q)$ is a complex Lie group
acting holomorphically and transitively on $Q$. The domain
$E$ is open and
dense in $Q$ and the set $Q\backslash E$ of the `points
at infinity' is an analytic subset of $Q$, but not a complex
submanifold in general.

The boundary $\partial D$
of $D$ is smooth only in the very special case, where also
$\|\!\cdot\!\|_\infty$ is a Hilbert
norm. Nevertheless, $\partial D$ is a finite union
of $\Gamma$-orbits, which are locally closed CR-submanifolds of $E$.
Every $K$-orbit $M$ in $\partial D$ is an Hermitian CR-submanifold of $E$
with respect to the metric induced from $E$, where $K$ acts by
CR-isometries.
We start with an orbit of a special nature: Denote by
$S=S(D)$ the set of all extreme points of the closed convex
set $\overline D$.
The following two statements are well known, but will
also be obvious from our discussion below.
\Lemma{SCHI} $S$ is a connected generic CR-submanifold of $E$.
Moreover, $S$ is the only compact $\Gamma$-orbit in $\overline D$,
consists of all $e\in\overline D$ with $K(e)=\Gamma(e)$ and
coincides with the Shilov boundary of $D$.
\qex
\Lemma{TUBE} The CR-submanifold $S$
is totally real if and only if
$D$ is biholomorphically equivalent to a `tube domain'
$\{z\in\CC^n:\Re(z)\in\Omega\}$, where $\Omega\subset\RR^n$
is an open convex cone. In this case $D$ is said to be of tube type.
\qex
\smallskip
A bounded symmetric domain $D$ is called {\sl irreducible} if it
is not biholomorphically equivalent to a direct product of
complex manifolds of lower dimensions. This is
known to be equivalent to
$K\subset\GL(E)$ acting irreducibly on $E$. There exists
(up to order) a unique representation of $D$ as
direct product $D=D_1\times\cdots\times D_k$, where all
$D_j$ are irreducible bounded symmetric domains
and are of the form $D_j=E_j\cap D$ for linear subspaces
$E_j\subset E$ with $E=E_1\oplus\cdots\oplus E_k$. Also, there exist
direct product representations $S(D)=S(D_1)\times\cdots\times S(D_k)$
for the Shilov boundaries and
$K=\GL(D_1)^0\times\cdots\times \GL(D_k)^0$. We call the
$D_j$ the {\sl factors} of $D$.
\smallskip

We are now able to formulate the main result of this section.
\Theorem{MINI} Let $D$ be a bounded symmetric domain.
Then the Shilov boundary $S$ of $D$ is a symmetric CR-manifold
and the following conditions are equivalent.
\0 The Levi cone of $S$ has non-empty interior at every point.
\1 $S$ is a minimal CR-manifold.
\1 Every smooth CR-function $f$ on $S$ has a unique holomorphic extension to $D$
        that has the same smoothness degree on $\overline D$ as $f$.
\1 $\Aut(D)=\AUT(S)$.
\1 $D$ does not have a factor of tube type.
\Formend
\noindent
For the proof \Ruf{ZWPR} we use the Jordan
theoretic approach to bounded symmetric domains as originated
by {\Icke Koecher} \Lit{KOEC}, for details in the following always compare
\Lit{LOSO}: There
exists a Jordan triple product $E^3\to E$,
$(x,y,z)\mapsto \{xyz\}$, that contains the full structural
information of $D$. This triple product is symmetric bilinear in the
outer variables $(x,z)$, conjugate linear in the inner variable $y$
and satisfies certain algebraic and spectral properties.
The group $\GL(D)$ of all linear $\|\!\cdot\!\|_\infty$-isometries of $E$ coincides
with the group of all linear triple automorphisms, more precisely
$$\GL(D)\;=\;\big\{g\in\GL(E):g\{xyz\}=\{gx\,gy\,gz\}
\steil{for all}x,y,z\in E\big\}\,.$$

An element $e\in E$ is called a {\sl tripotent} if $\{eee\}=e$
holds. The set $\Tri(E)$ of all tripotents in $E$ is a compact
real-analytic submanifold
of $E$ and the group $K$ acts transitively on
every connected component of $\Tri(E)$. Except for
$\{0\}$ every other connected component
of $\Tri(E)$ has positive dimension and is contained in
$\partial D$.
Tripotents may also be characterized geometrically
as `affine symmetry points' of $\overline D$ in the following
sense.
\Proposition{CHAR} The element $a\in\overline D$ is a tripotent if and only if there
exists an operator $\sigma\in\GL(D)$ with
\0 $\sigma(a)=a$,
\1 $\sigma(v)=-v$ for all $v\in E$ with
$\|a+tv\|\le1$ for all $t\in\TT$.
\Formend\noindent
We will postpone the proof of this criterion and
fix a tripotent $e\in E$ for a moment. The triple multiplication operator
$\mu=\mu_e\in\5L(E)$ defined by $z\mapsto \{eez\}$ is Hermitian
and splits $E$ into an orthogonal sum $E=E_1\oplus E_{1/2}\oplus E_0$
of eigenspaces to the eigenvalues $1,1/2,0$, called the {\sl Peirce spaces}
of the tripotent $e$. The canonical projection $P_k\colon E\to E_k$ maps
$D$ into itself and clearly
is a polynomial in $\mu$, more precisely
$$P_1=\mu(2\mu-1)\,,\qquad P_{1/2}=4\mu(1-\mu)\,,
\qquad P_0=(1-\mu)(1-2\mu)\,.\Leqno{PROJ}$$
The `Peirce reflection'
$\rho\colon=\exp(2\pi i\mu)=P_1-P_{1/2}+P_0$
is contained in $K$, fixes $e$ and leaves $\Tri(E)$ invariant.
In particular, also the projection $P_1+P_0$ maps $D$ into itself.

The tripotent $e\ne0$
is called {\sl minimal} if $E_1=\CC e$ holds and
is called {\sl maximal} if $E_0=0$ holds. For instance, the Shilov boundary $S$
of $D$ is just the set of all maximal tripotents.
$E$ becomes a complex {\sl Jordan algebra} (depending on the tripotent $e$)
with respect to the
commutative product $a\circ b\colon=\{aeb\}$, and $e^2\colon=e\circ e=e$
is an idempotent in $E$.
The Peirce space $E_1$ is a unital complex Jordan subalgebra
with identity element
$e$ and conjugate linear algebra involution $z\mapsto z^*\colon=\{eze\}$.
For every $a\in E_1$ and powers inductively
defined by $a^{k+1}\colon=a^k\circ a$, $\;a^0\colon=e$,
the linear subspace $\CC[a]\subset E_1$ is a
commutative, associative subalgebra (notice that the Jordan algebra
$E_1$ is not associative in general). The element $a\in E_1$ is called {\sl invertible}
if $a$ has an inverse $a^{-1}\in\CC[a]$.
The selfadjoint part $A\colon=\{z\in E_1:z^*=z\}$ of $E_1$ is a {\sl formally
real Jordan algebra}, i.e. a real Jordan algebra such that $x^2+y^2=0$
implies $x=y=0$ for all $x,y\in A$. Clearly, $E_1=A\oplus iA$ holds since the involution
is conjugate linear. For all $z\in E_1$ we denote by
$\Re(z)\colon=\Quot{(z+z^*)}{2}\in A$ the real part of $z$.
The set $Y\colon=\{a^2:a\in A\}$
of all squares in $A$ is a closed convex cone with $A=Y-Y$ and
$Y\cap-Y=\{0\}$. The interior 
$$\Omega\colon=\hbox{\rm Interior of }Y$$ coincides with $\exp(A)\subset A$
and also with the set of all $a\in Y$ that are invertible in $A$.
$\Omega$ is an open convex linearly-homogeneous cone in $A$.
The sesqui-linear
mapping $\Phi\colon E_{1/2}\oplus E_{1/2}\to E$ defined by
$\Phi(u,v)=2\{euv\}$ takes values in $E_1$ and satisfies
$\Phi(z,z)\in\overline\Omega$ for all $z\in E_{1/2}$, and
$\Phi(u,u)=0$ if and only if $u=0$.
To indicate the dependence on the tripotent $e\in E$ we also write
$E_k(e)$, $k=1,1/2,0$, for the Peirce spaces as well as
$A(e)$, $Y(e)$, $\Omega(e)$, $\rho_e$ and $\Phi_e$. Let us illustrate
these objects by a typical example.
\Example{FRUN} Fix arbitrary integers $p\ge q\ge1$ and consider the
complex Hilbert space $E\colon=\CC^{p\times q}$ of dimension
$n=pq$. Then
$D\colon=\{z\in E:\One-z^*\!z>0\}$ is a bounded symmetric
domain in $E$, where $\One=\One_q$ is the $q\times q$-unit matrix. $\|z\|_\infty^2$
is the largest eigenvalue of the Hermitian matrix $z^*\!z$, i.e. $\|z\|_\infty$
may be considered as the operator norm of $z$ if considered as operator
$\CC^q\to\CC^p$. The triple
product is given by $\{xyz\}=(xy^*\!z+zy^*\!x)/2$ and
$K\subset\GL(E)$ is the subgroup of all transformations $z\mapsto uzv$ with
$u\in\U(p)$ and $v\in\U(q)$. The Hilbert norm on $E$ given by
$\|z\|^2=\tr(z^*\!z)$ is $K$-invariant. $\Tri(E)$ is the disjoint union of the
$K$-orbits $S_0,S_1,\dots,S_q$, where $S_k$ is the set of all tripotents
$e\in E$ that have matrix rank $k$. In particular, if we write
every $z\in E$ as block matrix ${ab\choose cd}$ with
$a\in\CC^{k\times k}$ and matrices $b,c,d$ of suitable sizes,
then $e={\one_{\]4q}0\choose 0\;0}$ is a tripotent in $S_k$.
The corresponding Peirce spaces $E_1$, $E_{1/2}$ and $E_0$
consist of all matrices of the forms ${a0\choose00}$, ${0b\choose c0}$
and ${00\choose0d}$ respectively. Furthermore, $A$ is the real
subspace of all Hermitian matrices in $E_1$ and $\Omega\subset A$ is the convex
cone of all matrices ${a0\choose00}$ with $a\in\CC^{k\times k}$ positive
definite Hermitian. For every $u={0b\choose c0}\in H_eS_k$ we have
$\Phi_e(u,u)=2\{euu\}={a0\choose00}$ with $a=bb^*+c^*\!c$.
Finally, $S=S_q$ is the Shilov boundary of $D$. $S$ consists of all
matrices in $E$ whose column vectors are orthogonal in $\CC^p$, or
equivalently, which represent isometries
$\CC^q\to\CC^p$. The group $\Gamma$ is the set of all transformations
$$z\longmapsto(\alpha z+\beta)(\gamma z+\delta)^{-1}
\Steil{with}\pmatrix{\alpha&\beta\cr\gamma&\delta\cr}\in\SU(p,q)$$
and $\alpha,\beta,\gamma,\delta$ matrices of sizes $p\times p$,
$p\times q$, $q\times p$ and $q\times q$ respectively. In case $p=q>1$
the groups $\GL(D)$ and $\Aut(D)$ have two connected components, in all
other cases these groups are connected.
-- For the
special case $q=1$ we get for $D$ the euclidean ball $\BB$ in $E=\CC^p$
with Shilov boundary the unit sphere $S=S_1=\partial D$ as studied in
Example \ruf{EXAM}. For every $e\in S$ then $E_1(e)=\CC e$ holds
and $E_{1/2}(e)$ is the orthogonal complement of $e$ in the Hilbert space $E$.
\qex

\bigskip
Two tripotents $e,c\in E$ are called (triple) {\sl orthogonal} if $c\in E_0(e)$
holds. Then also $e\in E_0(c)$ is true and $e\pm c$ are
tripotents. An ordered tuple $(e_1,e_2,\dots,e_r)$ of pairwise orthogonal minimal
tripotents in $E$ is called a {\sl frame} in $E$ if there does not exist
a minimal tripotent $e\in E$ that is orthogonal to all $e_j$ in the triple
sense. All frames in $E$ have the same length $r$, which is called the
{\sl rank} of the bounded symmetric domain. Every element $a\in E$ has a representation
$$a=\lambda_1e_1+\lambda_2e_2+\cdots+\lambda_re_r,\qquad\qquad
\|a\|_\infty=\lambda_1\ge\lambda_2\ge\cdots\ge\lambda_r\ge0\,,\Leqno{SING}$$
where $(e_1,e_2,\dots,e_r)$ is a frame depending on $a$. The real numbers
$\lambda_j=\lambda_j(a)$ are uniquely determined by $a$ and are
called the {\sl singular values} of $a$. In general, the frame
$(e_1,e_2,\dots,e_r)$ is not uniquely determined by $a$. For every
$a\in\overline D$ there is a unique
representation $$a=e+u\Steil{with}e=\colon\epsilon(a)\in\Tri(E)
\steil{and}u\in D\cap E_0(e)\,.\Leqno{TYPI}$$
The Shilov-boundary of $D$ is given by
$$S=\{a\in E: \lambda_1(a)=\lambda_2(a)=\cdots=\lambda_r(a)=1\}\,.\Leqno{SHIL}$$
In case $D$ is irreducible,
the compact group $K$ acts transitively on the set of all
frames and hence any two elements $a,b\in E$ are in the same
$K$-orbit if and only if $\lambda_j(a)=\lambda_j(b)$ holds for all $j$.

These considerations can be used to prove the following property.

\Proposition{RATH} If $S$ is totally real, it is rationally convex.

\Proof For every $e\in S$, the Jordan algebra $E=E_1(e)$ has $e$ as the unit element. 
It is known (compare for instance \Lit{BRKO} or \Lit{MCCR}) that there exists a unique
homogeneous polynomial function $N\colon E\to\CC$ of degree $r$ 
such that the following is satisfied:
\0 $z\in E$ is invertible if and only if $N(z)\ne0$ \quad and \quad
(ii) $N(e)=1$.

\noindent $N$ is called the (generic) norm of the unital Jordan algebra $E$.
It is known that there exists a character $\chi\colon K\to\TT$ such that
$N(gz)=\chi(g)N(z)$ holds for all $g\in K$ and all $z\in E$. On the
other hand, for every frame $(e_1,\dots,e_r)$ in $E$ with $e_1+\cdots+e_r=e$
and every complex linear combination $z=z_1e_1+\cdots+z_re_r$ we have
$N(z)=z_1z_2\cdots z_r$. 
This implies the following characterization of the Shilov boundary in the tube type case.
$$S=\{z\in\overline D:|N(z)|=1\}\,.$$
In particular, for every $a\in\overline D\backslash S$,
the rational function $(N-N(a))^{-1}$ is holomorphic in a neighbourhood
of $S$ and has no holomorphic extension to $a$, i.e. the rational convex hull of $S$ in $E$ coincides with $S$.
\qed

\noindent The Shilov boundary $S$ of $D$ in Example \ruf{FRUN} is totally real
if and only if $p=q$ holds, and then $S=\U(q)$ is the unitary group. For the
unit matrix $e\in E=\CC^{q\times q}$ the Jordan product on $E$ is given by
$a\circ b=(ab+ba)/2$ and invertibility in the Jordan sense is the same as
in the associative sense. In particular, $N(z)=\det(z)$ is the norm of $E$.

\Joker{ERPR}{Proof of Proposition \CHAR} In case $a$ is a tripotent, every $v\in E$
with $\|a+tv\|\le1$ for all $t\in\TT$
is contained in $E_0(a)$ and we may take
$\sigma\colon=-\exp(\pi i\mu_a)=P_1-iP_{1/2}-P_0\in K$, where $\mu_a$ is
the triple multiplication operator $z\mapsto \{aaz\}$ on $E$.
Conversely, suppose that $a$ satisfies
\ruf{CHAR}.i-ii and write $a=e+u$ as in \ruf{TYPI}.
Then $\sigma(u)=-u$ follows from the assumptions.
For every $t>1$ with $tu\in\overline D$
we have $a-(1+t)u=e-tu\in\overline D$
and hence $\sigma(a-(1+t)u)=e+(2+t)u\in\overline D$, i.e $\;(t+2)u\in\overline D$
and hence $u=0$. Therefore, $a=e$ is a tripotent.
\qed
\medskip
Fix a frame $(e_1,e_2,\dots,e_r)$ in $E$ and consider for all
integers $0\le i,j\le r$, the {\sl refined Peirce spaces}:
$$E_{ij}\colon=\big\{z\in E\;:\;2\{e_ke_kz\}=(\delta_{ik}+
\delta_{kj})z\steil{for}1\le k\le r\big\}\,.$$
Then, if we put $e_0\colon=0$,
$$E=\bigoplus_{0\le i\le j\le r}\!\!E_{ij}\;,\qquad E_{ii}=\CC\[2e_i
\Steil{and}\{E_{ij}E_{jk}E_{kl}\}\subset E_{il}$$
hold for all $0\le i,j,k,l\le r$. Also, $\{E_{ij}E_{kl}E\}=0$
if the index sets $\{i,j\}$ and $\{k,l\}$ are disjoint. Furthermore, $D$
has no tube type factor if and only if $E_{i0}\ne 0$ for $1\le i\le r$.
To indicate the dependence of $E_{ij}$ on the given frame
we also write $E_{ij}(e_1,e_2,\dots,e_r)$.

Now consider a $\Gamma$-orbit $\Sigma\subset\overline D$. Then it is known
that there is a tripotent $e$ in $E$ with $\Sigma=\Gamma(e)$
and that $T_e\Sigma=iA\oplus E_{1/2}\oplus E_0$ is the tangent space at $e\in \Sigma$,
where the Peirce spaces refer to the tripotent $e$.
This implies that $\Sigma$ is a homogeneous generic locally-closed CR-submanifold
of $E$ with holomorphic tangent space $H_e\Sigma=E_{1/2}\oplus E_0$. The orbit
$M\colon=K(e)$ is a compact submanifold of $\Sigma$ with tangent space
$T_eM=iA\oplus E_{1/2}$ and holomorphic tangent space $H_eM=E_{1/2}$.
Furthermore, $M=\Sigma\cap\Tri(E)$ and $\epsilon\colon\Sigma\to M$
(compare \ruf{TYPI}) is a fibre bundle with
typical fibre $D\cap E_0$.
\Lemma{RUPP} Every connected component $M$ of $\Tri(E)$ is a symmetric CR-manifold.
\Proof Fix an arbitrary element $e\in M$. For the decomposition of $D$ into a direct
product $D_1\times\cdots\times D_k$
of irreducible factors we get a decomposition $e=e_1+\cdots+e_k$
with tripotents $e_j\in E_j$ and a decomposition $M=M_1\times\cdots\times M_k$
with $M_j=\GL(D_j)^0(e_j)$, that is, we may assume without loss of generality
that $D$ is irreducible. To begin with, suppose that $M$ is totally real, i.e. $E_{1/2}=0$.
Then, by irreducibility, also $E_0=0$ holds and $M=\exp(iA)$ is the
`generalized unit circle'
in $E_1=E$. Furthermore, $s_e(z)=z^*$ leaves $M$ invariant and hence is a symmetry
of $M$ at $e$, i.e. $M$ is symmetric in this case. Now suppose, that $M$
is not totally real, i.e. $H_eM=E_{1/2}\ne0$. Then the Peirce reflection
$\rho_e$ maps $M$ into itself and satisfies
$H_eM\subset\Fix(-\rho_e)$. Therefore, as soon as we know that
$M$ is a minimal CR-manifold we know that $\rho_e$ is a symmetry of
$M$ at $e$ and hence that $M$ is symmetric. For the minimality of $M$
it is enough to show that $H_e^2M=iA$ holds, where $H_e^2M$ is as in
section \ruf{zwei}. But this is a consequence of
the following Proposition \ruf{FOLL}.\qed

\Proposition{FOLL} Let $e$ be a tripotent in $E$ and denote by
$M$ the connected $e$-component of $\Tri(E)$. Then
$\,H_e^2M\subset iA(e)$ and the Levi form $E_{1/2}(e)\times E_{1/2}(e)\to E_1(e)$
of $M$ at $e$ is given by $(u,v)\mapsto-2\{euv\}$, i.e. $L-e=-\Phi_e$.
In case $D$ has no tube type factor,
the convex hull of $\{\Phi_e(u,u):u\in E_{1/2}(e)\}$
in $A(e)$ has the cone $-\Omega(e)$ as interior and then, in particular,
$H_e^2M=iA(e)$ holds.
\Proof For every $u\in H_eM=E_{1/2}(e)$ define the vector field $X^u$ on $E$ by
$X_a^u=4\{aau\}-4\{aa\{aau\}\}$ for all $a\in M$. Then $X_e^u=u$
and $X_a^u\in H_aM$ for all $a\in M$ by \Ruf{PROJ}. A simple calculation
gives $[X^u,X^v]_e=2\{evu\}-2\{euv\}\in iA(e)$. This shows that $-\Phi_e$ is the
Levi form at $e\in M$.
Let $C$ be the convex hull of $\{\Phi_e(u,u):u\in E_{1/2}(e)\}$.
Then $C\subset\overline\Omega(e)$ is clear. For the proof of
the opposite inclusion
fix an arbitrary element $a\in\overline\Omega(e)$.
Then there exists an integer $k\le r$ and a  representation
$a=\lambda_1e_1+\cdots+\lambda_ke_k$, where $(e_1,\dots,e_k)$ is
a family of pairwise orthogonal minimal idempotents
in the formally real Jordan algebra $A(e)$ summing up to $e$ and where
all coefficients $\lambda_j$ are $\ge0$. This means that we only
need to show that $e_j\in C$ for $1\le j\le k$. For this we
extend $(e_1,\dots,e_k)$ to a frame $(e_1,\dots,e_r)$ of $E$ and fix $j\le k$. Since
by assumption $D$ has no tube type factor we have $E_{j0}\ne0$. But then
$\Phi_e(u,v)=2\{e_juu\}$ cannot vanish for all $u,v\in E_{j0}$ since otherwise
there would exist a tripotent $c\ne0$ in $E_{j0}$ that is orthogonal
to all $e_i$, $1\le i\le r$. This implies $e_j\in C$.
\qed

\medskip
\Joker{ZWPR}{Proof of Theorem \MINI}
\Proof $S$ is symmetric by Lemma \ruf{RUPP} since $S$ is a connected component
of $\Tri(E)$. Fix an element $e\in S$. Then with $Q$ as in \ruf{COMP}
there exists an automorphism $\gamma\in\Aut(Q)$, called {\sl Cayley
transformation,} mapping $D$ biholomorphically onto the Siegel domain
$$\HH\colon=\big\{(t,v)\in E_1\oplus E_{1/2}:t+\overline t-\Phi(v,v)\in\Omega\big\}$$
in $E=E_1\oplus E_{1/2}$, where the Peirce spaces $E_k$, the cone $\Omega\subset A$
and the Hermitian map $\Phi\colon E_{1/2}\times E_{1/2}\to E_1$ refer to 
the tripotent $e$.
The transformation $\gamma$ satisfies $\gamma^4=\id$,\hfill\break
$S\cap\Fix(\gamma)=\{\pm ie\}$, $\gamma(-e)=0$ and is given by
$$\gamma(t,v)=\big((e-t)^{-1}\!\circ(e+t),\sqrt 2\,(e-t)^{-1}\!\circ v\big)\,,$$
where $(e-t)^{-1}$ is the inverse in the unital Jordan algebra $E_1$. The
domain $$V\colon=\big\{(t,v)\in S:(e-t)\;\steil{is invertible in}E_1\big\}$$ is dense in $S$
and $\gamma$ defines a CR-diffeomorphism from $V$ onto the CR-submanifold
$$N\colon=\big\{(t,v)\in E_1\oplus E_{1/2}:t+\overline t=\Phi(v,v)\big\}\;\subset\;\partial\HH$$
of $E$, compare also Example \ruf{EXAX}.
\hfill\break
\To12 is an immediate consequence of the definitions (see \ruf{LEVI}) and holds
for every CR-manifold.
\hfill\break
\To23. 
Suppose, $S$ is minimal. Since the Shilov boundary of a bounded
symmetric domain of tube type is totally real, 
$D$ cannot have a factor of tube type.
Then by Proposition \ruf{FOLL}, the interior of 
$\{\Phi(v,v):v\in E_{1/2}\}=\{L(v,v):v\in E_{1/2}\}$ coincides with the cone $\Omega$,
where $L$ denotes the Levi form of $N$ at $0$ with respect to the obvious identification
$(T_0N/H_0N)\otimes\CC \cong E_1$. Let $v\in\Omega$ be an arbitrary vector.
By the extension result of \Lit{BP}, every CR-function $f$ on $N$ extends holomorphically
to a small wedge in the direction $v$, in particular, to a neighbourhood of a subset of the type $(N+ \RR_+ v)\cap U$,
where $U$ is a neighbourhood of $0$ in $E_1\oplus E_{1/2}$.
Furthermore, the wedge extension is of the same smoothness degree as $f$
(due to \Lit{BCT}, see also \Lit{BER}, Theorem~7.5.1,
since $f$ is of slow growth by the Cauchy estimates).
Using the transformations $(t,v)\mapsto(st,s^2v)$, $s>0$,
we see that $f$ automatically extends to a neighbourhood of $(N+ \RR_+ v)$.
Since $v\in\Omega$ is arbitrary, $f$ extends holomorphically to the whole of $\HH$. 
This implies via the Cayley transformation that
every CR-function $f$ on $S$ has a smooth extension to $D\cap S$ which is holomorphic on $D$. 
It remains to prove that the extension of $f$ is of the same smoothness degree on the boundary $\partial D$.
For every $0<r<1$, define the continuous function $f_r$ on $\overline D$
by $f_r(z):=f(rz)$. Then for $r\to1$, the functions $f_r$ converge uniformly on $S$ to $f$.
Since $S$ is the Shilov boundary of $D$ the convergence is also uniform on $\overline D$,
i.e. $f$ extends continuously to $\overline D$. The smoothness is obtained
by the same argument applied to the partial derivatives.
\hfill\break
\To34. $\Aut(D)\subset\AUT(S)$ follows from \Ruf{COMP}. Assume (ii) and consider a
transformation $g\in\AUT(S)$. Then $g$ extends to a continuous mapping
$g\colon\overline D\to E$ which is holomorphic on $D$. As a consequence of the
maximum principle, $g(\overline D)$ is contained in the closed convex hull
of $g(S)=S$, which is $\overline D$. By the same argument, $h\colon=g^{-1}$ extends
to a continuous map $h\colon\overline D\to\overline D$ which is holomorphic on $D$.
Then $h\circ g=g\circ h=\id$ shows $g\in\Aut(D)$.
\hfill\break
\To45. Suppose $D$ has a factor of tube type. Then $S$ is a direct product of
a CR-manifold with a totally real CR-manifold of positive dimension. In particular,
$\AUT(S)$ cannot be a Lie group of finite dimension like $\Aut(D)$.
\hfill\break
\To51 follows from Proposition \ruf{FOLL}.
\qed

Theorem~\ruf{MINI} together with Proposition~\ruf{RATH} can be used to 
calculate both polynomial and rational convex hulls of $S$ explicitly.
In particular, they are finite unions of disjoint connected real-analytic CR-submanifolds
(forming a stratification in the sense of Whitney).
We call a smooth function on such a union a CR-function if it is CR on each single CR-submanifold
(this notion is independent of the partition into CR-submanifolds).

\Corollary{HULLS} Let $E=E_1\oplus E_2$ be the canonical splitting such
that $D_1:=D\cap E_1$ is of tube type and $D_2:=D\cap E_2$ has no tube type factor.
Denote by $S_1\subset \partial D_1$ and $S_2\subset \partial D_2$ the corresponding Shilov boundaries.
Then the following holds.
\0 Both convex and polynomial convex hulls of $S$ coincide with $\overline D$.
\1 The rational convex hull $\hat S$ of $S$ is given by $\hat S = S_1\times\overline D_2$.
\1 Every smooth CR-function $f$ on $S$ extends uniquely to a CR-function on $\hat S$ of the same smoothness degree.

\Proof Since $S$ is the Shilov boundary of $D$, 
$|P(z)|\le\|P\|_S$ holds for every holomorphic polynomial $P$ and every $z\in \overline D$. 
Hence $\overline D$ is contained in the polynomial convex hull of $S$.
The latter is always contained in the convex hull of $S$, which is $\overline D$. 
This proves (i) (the statement about the convex hull also follows from the classical Krein-Milman theorem).

For the rational convex hull, we obtain $\hat S_1 = S_1$ by Proposition~\ruf{RATH}.
This shows $\hat S\subset S_1\times\overline D_2$.
On the other hand, every rational function on $E_2$, holomorphic in a neighbourhood of $S_2$,
is continuous on $\overline D_2$ by Theorem~\ruf{MINI}. This implies the opposite inclusion 
$\hat S\supset S_1\times\overline D_2$ and therefore (ii).

Finally, let $f$ be a CR-function on $S$.
Then, for every $z_1\in S_1$, Theorem~\ruf{MINI} guarantees that 
$f$ has a unique smooth extension $\hat f$ to $\{z_1\}\times\overline D_2$ which is holomorphic on $\{z_1\}\times D_2$.
By the smoothness, $\hat f$ is CR on each CR-submanifold of the boundary $\{z_1\}\times\partial D_2$.
To prove the smoothness of $\hat f$ on $\hat S$, we fix a convergent sequence $z_1^m\to z_1^0$.
Then $\hat f(z_1^m,\cdot)$ converges to $\hat f(z_1^0,\cdot)$ uniformly on $S_2$ and therefore on $\partial D_2$,
because $S_2$ is the Shilov boundary. This shows that $\hat f$ is continuous.
The same argument applied to the partial derivatives of $\hat f$ shows that $\hat f$ is of the same
smoothness degree as $f$. Since $S_1$ is totally real, 
the holomorphic tangent spaces to every CR-submanifold of $\hat S$ are contained in $E_2$.
This shows that $\hat f$ is CR and finishes the proof of (iii).
\qed

\vfill\eject     %%%  #######  %%%

\bigskip
From the classification of all irreducible bounded symmetric
domains into the 6 types {\bf I, II,\dots,VI}
(compare f.i. \Lit{LOSO} \p4.11) it follows that there are precisely
the following 3 types of irreducible non-tube domains:
\medskip
{\Bf I$_{q,p}$} with $p>q\ge1$ arbitrary integers.
Then, as in Example \ruf{FRUN}, $E=\CC^{p\times q}$ and
$D=\{z\in E:\One-z^*\!z>0\}$ is the bounded symmetric
domain of rank $r=q$, where $\One$ is the $q\times q$-unit matrix.
The Shilov boundary of $D$
is the set $S\colon=S_q$ of all matrices in $E$
whose column vectors are orthogonal in $\CC^p$, i.e.
$$S=\{z\in\CC^{p\times q}:z^*\!z=\One\}\,.\Leqno{ZITA}$$
On $S$ the group $\SU(p)$
acts transitively by matrix multiplication from the left
with isotropy subgroup $\One\times\SU(p\!-\!q)$ at
$e\colon={\one\choose0}\in S$, i.e. $S=\SU(p)/(\One\times\SU(p\!-\!q))$ is
simply connected, has CR-dimension $(p-q)q$ and CR-codimension
$q^2$. Also, $\AUT(S)\approx\U(p,q)/\TT$ is connected.
Every closed subgroup $L\subset\U(q)$ acts freely on $S$ by matrix
multiplication from the right and $S/L$ again is a symmetric CR-manifold
of the same CR-dimension
in a natural way. For $L=U(q)$ we get the Grassmannian of
all $q$-planes in $\CC^p$ which is the reduction of $S$
as defined in section \ruf{fuenf}. The typical fibre of
the reduction map is the group $\U(q)$.
\qex\medskip
{\Bf II$_p$} with $p=2q+1$ an arbitrary odd integer $>3$. Let
$E\colon=\{z\in\CC^{p\times p}:z'=-z\}$ and define the bounded
domain $D\subset E$ as well as the Jordan triple product
by the same formulae as for {\bf I}$_{p,q}$. Then again $D$ is a bounded
symmetric domain of rank $r=q$ and
$\Gamma\subset\GL(E)$ is the group of all
transformations $z\mapsto uzu'$ with $u\in\U(p)$. For
$j\colon={\;0\hskip6pt\one\choose -\one\;0}\in\CC^{2q\times 2q}$
the matrix $e\colon={j\;0\choose0\;0}\in\CC^{p\times p}$ is in $S\colon=S_q$
and $\Gamma$ has isotropy subgroup $\Sp(q)\times\TT$ at $e$. Therefore
$S=SU(p)/(\Sp(q)\times1)$ is simply connected. The holomorphic tangent space
$H_eS$ is the space of all $z={0\;u\choose v\;0}\in E$ with $u\!=\!-v'\in\CC^{2q}$,
i.e. $S$ has CR-dimension $2q$ and CR-codimension $q(2q-1)$.
The reduction is the projective space $\PP_{2q}(\CC)=\SU(p)/S(\U(2q)\times\TT)$
and $\SU(2q)/\Sp(q)$ is the typical reduction fibre.
The group $\AUT(S)\;\approx\;\SO^*(2p)/\{\pm1\}$ is connected (compare
\Lit{HELG} \p451 and \p518 for the non-compact type D III).
\qex\medskip
{\Bf V} Here $D$ is the exceptional bounded symmetric domain of dimension
16 (non-compact type E~III on \p 518 of \Lit{HELG}). $D$ has rank 2 and
the Shilov boundary $S\colon=S_2$ has CR-dimension 8 and CR-codimension 8.
On $S$ the group $\Spin(10)$ acts transitively and the reduction $\widetilde S$ of $S$
is the symmetric Hermitian manifold $\SO(10)/(\SO(2)\times\SO(8))$, the complex
nonsingular quadric of dimension 8.
The group $\AUT(S)$ is a non-compact simple exceptional real
Lie group of type $E_6$ and has dimension 78.
\qex

\medskip
As a generalization of Example \ruf{EXAA} also the dual of \Ruf{ZITA} (compare
\Ruf{DUAL} and section \ruf{fuenfa})
can be described explicitely. Fix $e={\one\choose0}\in S$
and denote by $\rho=\rho_e$ the corresponding Peirce reflection
of $E=\CC^{p\times q}$. Then $\Fix(\rho)=\CC^{q\times q}$ and
$\Fix(-\rho)=\CC^{(p-q)\times q}$. On
$$R=\{z\in\CC^{p\times q}:\rho\Kl z ^*\!z=\One\}$$
the group $\U(q,p\!-\!q)$ acts transitively from the left with
compact isotropy subgroup $\One\times\U(p\!-\!q)$ at
$e$. Therefore there is a unique $\U(q,p\!-\!q)$-invariant
Riemannian metric on $R$ which coincides on $T_eR$
with the one induced from $E$. The restriction of $\rho$
to $R$ is a symmetry of $R$ at $e$, i.e. is a symmetric CR-manifold. Again, every closed
subgroup $L\subset\U(q)$ acts freely on $R$
from the right and $S/L$ is a symmetric CR-manifold
of the same CR-dimension.
For $L=U(q)$ we get the bounded symmetric domain of type {\bf I}$_{q,p-q}$,
the reduction of $R$.

\eject

\bigskip\klein\parindent0pt
{\gross References}

\Ref{ALO}Alexander, H.: Polynomial approximation and hulls in sets of finite linear measure in $\CC^n$. {\sl Amer. J. Math.} {\bf 93}, 65--74 (1971).
\Ref{ALN}Alexander, H.: The polynomial hull of a rectifiable curve in $\CC^n$. {\sl Amer. J. Math.} {\bf 110} (4), 629--640 (1988).
\Ref{ANDR}Andreotti, A.; Fredricks, G.A.: Embeddability of real analytic Cauchy-Riemann manifolds. {\sl Ann. Scuola Norm. Sup. Pisa Cl. Sci.}, Serie IV, {\bf 6} (2), 285--304 (1979).
\Ref{BCT}Baouendi, M.S.; Chang, C.H.; Treves, F.: Microlocal hypo-analyticity and extension of CR functions. {\sl J. Differential Geom.} {\bf 18}, 331--391 (1983).
\Ref{BER}Baouendi, M.S.; Ebenfelt, P.; Rothschild, L.P.: {\sl Real Submanifolds in Complex Spaces and Their Mappings}. Princeton Math. Series, Princeton Univ. Press, 1998.
\Ref{BRMIN}Baouendi, M.S.; Rothschild, L.P.: Cauchy-Riemann functions on manifolds of higher codimension in complex space. {\sl Invent. Math.} {\bf 101} (1), 45--56 (1990).
\Ref{BG}Bloom,~T.; Graham,~I.: On type conditions for generic real submanifolds of $\CC^n$. {\sl Invent.~Math.} {\bf 40}, 217--243 (1977).
\Ref{BISH}Bishop, E. Holomorphic completions, analytic continuation, and the interpolation of semi-norms. {\sl Ann. of Math.}, II. Ser. {\bf 78}, 468--500 (1963).
\Ref{BOGG}Boggess,~A.: {\sl CR Manifolds and the Tangential Cauchy-Riemann Complex}. Studies in Advanced Mathematics. CRC Press.~ Boca Raton~ Ann Arbor~ Boston~ London 1991.
\Ref{BP}Boggess,~A.; Polking,~J.C.: Holomorphic extension of CR functions. {\sl Duke Math.~J.} {\bf 49}, 757--784 (1982). % extensions in the direction of the Levi form
\Ref{BRKO}Braun, H.; Koecher,  M.: {\sl Jordan-Algebren.} 
Berlin-Heidelberg-New York Springer 1966.
\Ref{CHIR}Chirka, E.M.: {\sl Introduction to the geometry of CR-manifolds}. Russ. Math. Surv. {\bf 46} (1), 95--197 (1991); translation from Usp. Mat. Nauk {\bf 46} (1), 81--164 (1991).
\Ref{DINH}Dinh,~T.C.: Enveloppe polynomiale d'un compact de longueur finie et chaines holomorphes \`a bord rectifiable. {\sl Acta Math.} {\bf 180} (1), 31--67 (1998).
\Ref{HELG}Helgason, S.: {\sl Differential Geometry, Lie Groups and Symmetric Spaces.} Academic Press.~ New York~ San Francisco~ London 1978.
\Ref{JACO}Jacobowitz, H.: {\sl An introduction to CR structures}. Mathematical Surveys and Monographs {\bf 32}. Providence, Rhode Island (1990).
\Ref{KANT}Kantor, I.L.; Sirota, A.I.; Solodovnikov, A.S.: Bisymmetric Riemannian spaces. Isv. Russ. Akad. Nauk Ser. Mat. {\bf 59}, 85--92 (1995); Russian Sci. Isv. Math. {\bf 59}, 963--970 (1995).
\Ref{KOEC}Koecher, M.: {\sl An elementary approach to bounded symmetric domains.} Rice Univ. 1969.
\Ref{LAW}Lawrence, M.G.: Polynomial hulls of rectifiable curves. {\sl Amer. J. Math.} {\bf 117} (2), 405--417 (1995).
\Ref{LSSS}Loos, O.: Spiegelungsr\"aume und homogene symmetrische R\"aume. {\sl Math. Z.} {\bf 99}, 141--170 (1967).
\Ref{LOSS}Loos, O.: {\sl Symmetric Spaces I/II.} W. A. Benjamin, Inc. ~New York~  Amsterdam 1969.
\Ref{LOSO}Loos, O.: {\sl Bounded symmetric domains and Jordan pairs.} Mathematical Lectures. Irvine: University of California at Irvine 1977.
\Ref{MCCR}McCrimmon, K.: The generic norm of an isotope of a Jordan algebra. Scripta Math. {\bf 29}, 229--236 (1973).
\Ref{STOE}Stolzenberg, G.: Polynomially and rationally convex sets. {\sl Acta Math.} {\bf 109}, 259--289 (1963). 
\Ref{STO}Stolzenberg, G.: Uniform approximation on smooth curves. {\sl Acta Math.} {\bf 115}, 185--198 (1966).
\Ref{TUCR}Tumanov, A.E.: {\sl The geometry of CR-manifolds}. Several complex variables. III. Geometric function theory, Encyclopedia Math. Sci. {\bf 9}, 201--221 (1989); translation from Itogi Nauki Tekh., Ser. Sovrem. Probl. Mat., Fundam. Napravleniya {\bf 9}.
\Ref{TUMA}Tumanov, A.E. Extension of CR functions into a wedge from a manifold of finite type. {\sl Math. USSR, Sb.} {\bf 64} (1), 129--140 (1989); translation from {\sl Mat. Sb., Nov. Ser.} {\bf 136(178)}, No.1(5), 128--139 (1988).
\Ref{TH}Tumanov,~A.E.; Khenkin,~G.M.: Local characterization of holomorphic automorphisms of Siegel domains. {\sl Funct. Anal. Appl.} {\bf 17}, 285--294 (1983); translation from Funkts. Anal. Prilozh. {\bf 17} (4), 49--61 (1983). % local CR homeom betwe
\Ref{WERA}Wermer, J.: The hull of a curve in $\CC^n$. {\sl Ann. of Math.}, II. Ser. {\bf 68}, 550--561 (1958).
\Ref{WERE}Wermer, J.: Polynomially convex hulls and analyticity. {\sl Ark. Mat.} {\bf 20}, 129--135 (1982).

\closeout\aux\bye